\documentclass[12pt]{amsart}
\usepackage{amsfonts,amsmath,amssymb,mathrsfs,verbatim,amsthm,amscd}
\newtheorem {thm}{Theorem}

\newtheorem {cor}[thm]{Corollary}

\newcommand{\mit}{}
\newcommand{\cal}{\mathcal}

\newcommand{\N}{\mathbb N}

\newcommand{\C}{\mathbb C}
\newcommand{\R}{\mathbb R}
\newcommand{\Z}{\mathbb Z}

\newcommand{\T}{\mathbb T}
\renewcommand{\L}{\mathbb L}
\renewcommand{\P}{\mathcal P}

\begin{document}

\title[Elliptic hypergeometric functions]{
\bf Classical elliptic hypergeometric \\
functions and their applications}

\thanks{
This work is supported in part by the Russian Foundation for Basic Research
(RFBR) grant no. 03-01-00780. \\  \indent
Proceedings of the International
Workshop {\em ``Elliptic integrable systems"}
(RIMS, Kyoto, November 8-11, 2004), Rokko Lectures in Mathematics Vol. 18
(Edited by M. Noumi and K. Takasaki),
Department of Mathematics, Kobe University, 2005, pp. 253--287.
(In the present version, we removed some missprints noticed in
the published paper and updated the references.) }

\author{V. P.  Spiridonov} 
\address{Max--Planck--Institut f\"ur Mathematik,
Vivatsgas\-se 7, D-53111, Bonn, Germany
(on leave from Bogoliubov Laboratory of Theoretical Phy\-sics,
JINR, Dubna, Moscow region 141980, Russia);
e-mail address: {\tt spiridon@theor.jinr.ru} }

\maketitle

\begin{flushright}
\em To the memory of A. A. Bolibrukh
\end{flushright}

\begin{abstract}

General theory of elliptic hypergeometric series and integrals is outlined.
Main attention is paid to the examples obeying properties
of the ``classical" special functions. In particular, an elliptic
analogue of the Gauss hypergeometric function and some of its
properties are described. Present review is based on author's
habilitation thesis \cite{spi:thesis} containing a more detailed
account of the subject.
\end{abstract}

\tableofcontents  

\section{General definition of univariate elliptic hypergeometric series
and integrals}

\noindent
\underline{B{\scriptsize ROAD DEFINITION}} ($n=1$, univariate case)
\cite{spi:theta1,spi:theta2}.
\par
Formal contour integrals $\int_C\Delta(u)du$ and series
$\sum_{n\in {\mathbb Z}}c_n$ are called elliptic hypergeometric integrals and series, if
there exist three constants $\omega_1,\omega_2,\omega_3\in\C$ such that
\begin{eqnarray*}&\bullet&
\Delta(u+\omega_1)=h(u)\Delta(u),\\
&& \hbox{where $h(u)$ is an elliptic function of $u$ with periods }
\omega_2,\, \omega_3, \mbox{ i.e.,}\\
&& \mbox{$h(u)$ is meromorphic and }\\
&& \qquad h(u+\omega_2)=h(u+\omega_3)=h(u),\quad \mbox{Im}(\omega_2/\omega_3)\ne 0;\\
&\bullet& c_{n+1}=h(n\omega_1)c_n, \\
&&\hbox{where $h(n\omega_1)$ is an elliptic function of $n$ with periods }
\frac{\omega_2}{\omega_1},\
\frac{\omega_3}{\omega_1}.
\end{eqnarray*}
\par
There is a functional freedom in the definition of integrals:
$\Delta(u)\to \varphi(u)\Delta(u)$,
where  $\varphi(u)$ is an arbitrary $\omega_1$-periodic function,
$\varphi(u+\omega_1)=\varphi(u)$ (such a freedom is not essential for
series).

\vspace{2mm}
\par\noindent
\underline{N{\scriptsize ARROW DEFINITION OF INTEGRALS}.}
\par
Formal contour integrals $\int_C\Delta(u)du$
are called elliptic hypergeometric integrals,
if $\Delta(u)$ is a meromorphic solution of three linear
first order difference equations
$$
\Delta(u+\omega_i)=h_i(u)\Delta(u),\quad i=1,2,3,
$$
where $h_i(u)$ are elliptic functions with the periods
$\omega_{i+1},\, \omega_{i+2}$ (we set $\omega_{i+3}=\omega_i).$

If all $h_i(u)\neq const$, then $\mbox{Im}(\omega_i/\omega_j)\neq 0,\,
i\neq j$. Interesting situations occur when one $h_i(u)=const$,
in which case we can have either $\mbox{Im}(\omega_i/\omega_{i+1})=0$ or
$\mbox{Im}(\omega_i/\omega_{i+2})=0$.
For pairwise incommensurate $\omega_i$, the functional freedom
in the definition of $\Delta(u)$ is absent due to the non-existence
of triply periodic functions.

Thus, we have in general three elliptic curves, but only
two of them are independent. One can consider also elliptic hypergeometric
functions in a more general context, when $h_i(u)$ are $N\times N$ matrices
with elliptic function entries.

It is possible to abandon the requirement of double periodicity of
$h(u)$ in favor of its double quasiperiodicity similar to the Jacobi
theta or Weierstrass sigma functions. This leads to
a more general family of theta hypergeometric series and integrals
(theta analogs of the Meijer function)
\cite{spi:theta1,spi:theta2}, but we skip their consideration.
In the next section we describe certain ``classical'' special
functions of hypergeometric type and their elliptic generalizations.

\section{An overview of classical hypergeometric functions}

The Euler's beta integral \cite{aar:special}
$$
\int_0^1 x^{\alpha-1}(1-x)^{\beta-1}dx
=\frac{\mit\Gamma(\alpha)\mit\Gamma(\beta)
}{\mit\Gamma(\alpha+\beta)},\quad {\rm Re}\ \alpha,\ {\rm Re}\ \beta >0,
$$
determines: i) the measure for Jacobi polynomials; ii) an
integral representation for the ${}_2F_1$ series.
Namely, Jacobi polynomials
$$
P_n(x)=\frac{(\alpha)_n}{n!} {}_2F_1
\bigg( {{-n,\, n+\alpha+\beta-1}\atop \alpha};x \bigg),
$$
where $(\alpha)_n=\alpha(\alpha+1)\cdots(\alpha+n-1)$ is the
Pochhammer symbol, satisfy the orthogonality relations
\begin{eqnarray*}
&& \langle P_n,P_m\rangle=\int_0^1x^{\alpha-1}(1-x)^{\beta-1}P_n(x)P_m(x)dx
\\ &&\makebox[4em]{}
=\frac{\delta_{nm}}{2n+\alpha+\beta-1}\,
\frac{\mit\Gamma(n+\alpha)\mit\Gamma(n+\beta)}{\mit\Gamma(n+\alpha+\beta-1)n!}.
\end{eqnarray*}
The Gauss hypergeometric function has the form
\begin{eqnarray*}
&& {}_2F_1 \bigg( {{a,b}\atop c};x
\bigg)=\sum_{n=0}^\infty \frac{(a)_n(b)_n}{n!(c)_n}x^n
\\ &&\makebox[2em]{}
=\frac{\mit\Gamma(c)}{\mit\Gamma(b)\mit\Gamma(b-c)}
\int_0^1 t^{b-1}(1-t)^{b-c-1}(1-x t)^{-a}dt,
\end{eqnarray*}
where we skip for brevity relevant constraints upon the parameters.
It defines a solution of the hypergeometric equation
$$
y''(x)+\bigg(\,\frac{c}{x}+\frac{a+b-c+1}{x-1}\bigg)y'(x)+
\frac{ab}{x(x-1)}y(x)=0,
$$
which is analytical near the origin $x=0$.

Two integrals described above fit into the general pattern $\int_C\Delta(x)dx$
with the kernel $\Delta(x)=\prod_{j=1}^k(x-x_j)^{\alpha_j}$
and some free parameters $x_j$ and $\alpha_j$. This kernel is characterized by the
condition that its logarithmic derivative $\Delta'(x)/\Delta(x)=R(x)$ is a rational
function of $x$. A very natural generalization of this criterion consists
in the requirement that the kernel $\Delta(x)$ satisfies the first
order linear finite difference equation
$\Delta(x+\omega_1)=R(x)\Delta(x)$ with rational $R(x)$ (such a treatment is
valid already for the $_2F_1$ function via the Mellin-Barnes integral
representation). By definition, we obtain general {\em plain
hypergeometric integrals} for which
$$
\Delta(x)= \frac {\prod_{j=0}^{s-1} \mit\Gamma(x/\omega_1+u_j)}
{\prod_{j=0}^{r}  \mit\Gamma(x/\omega_1+v_j)}\,
 \varphi(x)\, y^{x/\omega_1},\quad \varphi(x+\omega_1)=\varphi(x),
$$
and $\Gamma(x)$ is the Euler's gamma function.

The Pochhammer series have the form $\sum_{n=0}^\infty c_n$ with
$$
\frac{c_{n+1}}{c_n}=R(n)
=\frac{\prod_{j=0}^{s-1}(n+u_j)}{(n+1)\prod_{j=1}^{r}(n+v_j)}\,y,
$$
which leads automatically to the expression
$$
\sum_{n=0}^\infty c_n={}_sF_r \bigg( {{u_0,\ldots  ,u_{s-1}}\atop {v_1,\ldots  ,v_r}};y
\bigg)=\sum_{n=0}^\infty\frac{(u_0)_n\cdots(u_{s-1})_n}{n!(v_1)_n\cdots(v_r)_n}\,y^n.
$$
\noindent
These series admit confluence limits like
${}_s F_r(y)\sim {}_{s-1} F_r(u_0 y)$ for $u_0\to \infty$.
Their $q$-generalization has by definition $ c_{n+1}/c_n=R(q^n)$ for $q\in\mathbb{C}$
and some rational $R(x)$, which leads uni\-que\-ly to the series
$$
{}_s\varphi_r \bigg( {{t_0,\ldots  ,t_{s-1}}\atop {w_1,\ldots  ,w_r}};q;y
\bigg)=\sum_{n=0}^\infty\frac{(t_0;q)_n\cdots(t_{s-1};q)_n}
{(q;q)_n(w_1;q)_n\cdots(w_r;q)_n}\,y^n,
$$
where $(t;q)_n=(1-t)(1-tq)\cdots(1-tq^{n-1})$ denotes the $q$-Poch\-ham\-mer
symbol. This definition differs from the one given in
\cite{aar:special,gas-rah:basic} by the inversion $q\to q^{-1}$ and
appropriate change of notation for parameters.
For $t_i=q^{u_i}$, $w_i=q^{v_i}$ and $q\to 1^-$,
we formally have ${}_s\varphi_r (y)\to {}_sF_r(\tilde{y})$
for some renormalized value of the argument $\tilde y$.
In a similar way one can reconstruct the bilateral series $_sH_r$ and $_s\psi_r$.

\underline{Elliptic hypergeometric series},
directly derived from the definition given
in the first section, have the form (the unilateral case)
$$
{}_{r+1}E_r\bigg({t_0,t_1,\ldots  ,t_r\atop w_1,\ldots  ,w_r};q,p;y\bigg)
=\sum_{n=0}^\infty\frac{(t_0)_n(t_1)_n\cdots (t_r)_n}{(w_0)_n(w_1)_n
\cdots (w_r)_n}y^n,
$$
where $w_0=q$ (the canonical normalization) and
\begin{itemize}
\item [$\bullet$]
$(t)_n=\theta(t,tq,\ldots,tq^{n-1};p)\equiv\prod_{j=0}^{n-1}\theta(tq^j;p),$
\par \vspace{1mm}
$\qquad$ \framebox{
$ \theta(t;p)=(t;p)_\infty(pt^{-1};p)_\infty,
$} \par \vspace{1mm}
$\qquad
(t;p)_\infty=\prod_{n=0}^\infty(1-tp^n),\quad |p|<1;
$

\smallskip
\item [$\bullet$]
\framebox{$
\prod_{j=0}^r t_j=\prod_{j=0}^r w_j,
$} the balancing condition.
\end{itemize}
The elliptic Pochhammer symbol $(t)_n$ (denoted also in some other places
as $(t;q,p)_n$, $\theta(t)_n,$ or $\theta(t;p;q)_n$) degenerates to $(t;q)_n$ for
$p\to 0$, $(t)_n\to(t;q)_n$. Therefore for generic fixed
$t_j, w_j$ we have the termwise limiting relation
$_{r+1}E_r\to {}_{r+1}\varphi_r$ with the balancing restriction
indicated above (which does not coincide with the balancing condition
usually accepted for $q$-hypergeometric series \cite{gas-rah:basic}).

For $p=e^{2\pi i \tau},\, \mbox{Im}(\tau)>0,$ and any
$\sigma, u\in\mathbb{C}$, $q=e^{2\pi i\sigma}$,
we have the following relation between $\theta(t;p)$ and
the Jacobi $\theta_1(x)\equiv\theta_1(x|\tau)$
function
\begin{eqnarray*}
&&\theta_1(\sigma u| \tau)=
-i\sum_{k=-\infty}^\infty (-1)^kp^{(2k+1)^2/8}q^{(k+1/2)u}
\\ && \makebox[4em]{}
= ip^{1/8} q^{-u/2}\: (p;p)_\infty\: \theta(q^{u};p).
\end{eqnarray*}
Properties $\theta_1(x+1)=-\theta_1(x)$,
$\theta_1(x+\tau) = -e^{-\pi i\tau-2\pi ix}\theta_1(x)$ and
$\theta_1(-x)=-\theta_1(x)$ simplify to $\theta(pz;p)=\theta(z^{-1};p)
=-z^{-1}\theta(z;p)$.

For $_{r+1}E_r$ series we have
$$
\frac{c_{n+1}}{c_n}
=y\prod_{j=0}^r \frac{\theta(t_jq^n;p)}{\theta(w_jq^n;p)}
=h(n\omega_1),
$$
an elliptic function of $n\in\mathbb{C}$ with periods
$\omega_2/\omega_1$, $\omega_3/\omega_1$  for
\par\vskip 1mm
\centerline{
\fbox{$
q=e^{2\pi i\omega_1/\omega_2}, \qquad
p=e^{2\pi i\omega_3/\omega_2}. $} }
\smallskip
\noindent
The integer $r+1$ is called the order of $h(x)$ and it counts the number
of zeros or poles of $h(x)$ inside the fundamental parallelogram of periods.
Due to the balancing condition, we have an interesting (and useful) property
$$
{}_{r+1}E_r\bigg({t_0,t_1,\ldots  ,t_r\atop w_1,\ldots  ,w_r};q,p;y\bigg)
={}_{r+1}E_r\bigg({t_0^{-1},t_1^{-1},\ldots  ,t_r^{-1}\atop
w_1^{-1},\ldots  ,w_r^{-1}};q^{-1},p;y\bigg).
$$

In the table below we describe known special functions with properties
generalizing in a natural way the $_2F_1$ hypergeometric function
features. It is rather sketchy
and does not pretend on completeness. The $_2F_1$ series is the classical
special function investigated by such giants as Euler, Gauss, Jacobi,
Riemann and many other mathematicians. Its $q$-analogue has been proposed
by Heine as far back as 1850. However, until the relatively recent time
(landmarked by the appearance of quantum algebras from exactly solvable
models of statistical mechanics) $q$-special functions did not attract
much attention.

Chebyshev put forward general theory of orthogonal polynomials
which played a major role in the search of classical special functions.
Jacobi polynomials satisfy a three term recurrence relation and the
hypergeometric equation. The general discrete set of $_3F_2$ polynomials was
constructed by Chebyshev (I am indebted to R. Askey and A. Zhedanov for pointing
this fact to me). Their continuous analogues and the $_3\phi_2$ series generalizations,
known as Hahn polynomials, were proposed much later.
These polynomials satisfy a second order difference equation
(instead of the differential equation)  in their argument
lying on some non-trivial ``grids". The next level of
generalization is given by the Racah and Wilson polynomials described by special
$_4F_3$ series. In 1985, Askey and Wilson have found \cite{ask-wil:some}
the most general set of orthogonal polynomials with the self-duality
property. They are

\medskip
\begin{center}
\underline{\bf \scriptsize CLASSICAL SPECIAL FUNCTIONS OF HYPERGEOMETRIC TYPE}
\end{center}
$$\!\!\!\!\!\!
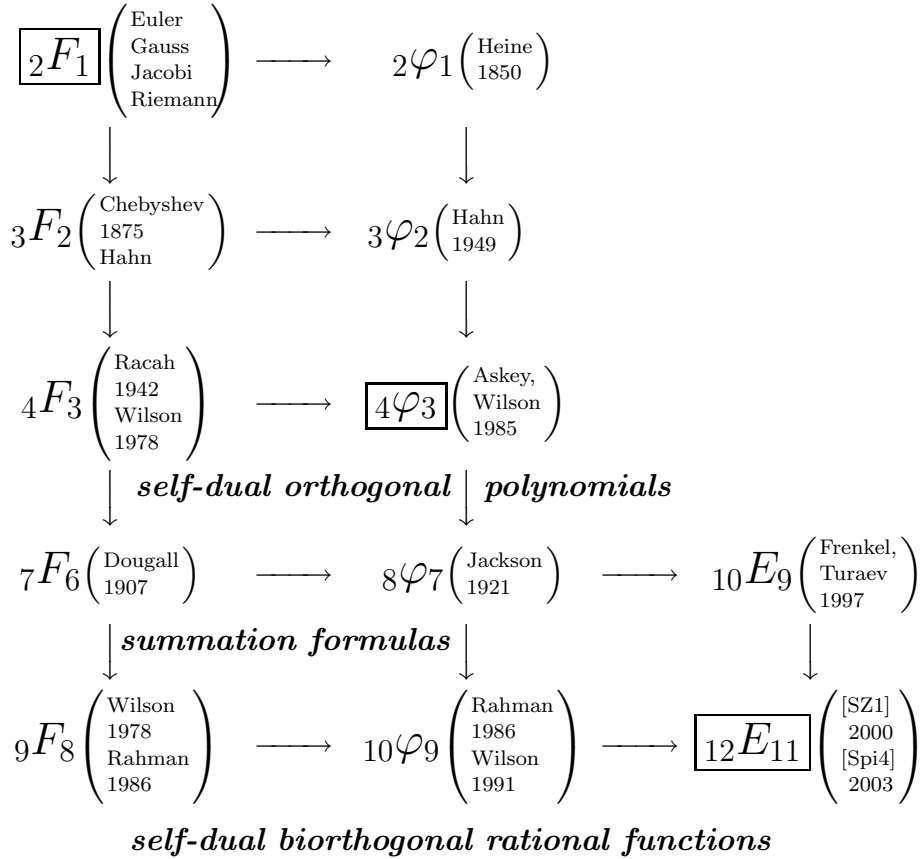
\begin{CD}
\text{\hskip 5mm}
\fbox{{\Large${}_2F_1$}}
{\tiny\left(\!\!
\begin{array}{l}
\text{Euler}\\
\text{Gauss} \\
\text{Jacobi}\\ 
\text{Riemann}
\end{array}\!\!\!\!
\right)}
@>>>
{\Large\text{
${}_2\varphi_1$}}
{\tiny\left(\!\!
\begin{array}{l}
\text{Heine}\\
\text{1850}
\end{array}\!\!
\right)}
@.\\
@VVV @VVV@.\\
{\Large\text{
${}_3F_2$}}
{\tiny\left(\!\!
\begin{array}{l}
\text{Chebyshev}\\
\text{1875}\\
\text{Hahn}
\end{array}\!\!
\right)}
\text{\hskip -1mm} 
@>>>
{\Large\text{
${}_3\varphi_2$}}
{\tiny\left(\!\!
\begin{array}{l}
\text{Hahn}\\
\text{1949}
\end{array}\!\!
\right)}
\text{\hskip 7mm}
@.\\
@VVV @VVV@.\\
{\Large\text{
${}_4F_3$}}
{\tiny\left(\!\!
\begin{array}{l}
\text{Racah}\\
\text{1942}\\
\text{Wilson}\\
\text{1978}
\end{array}\!\!
\right)}
@>>>
{\Large
\fbox{${}_4\varphi_3$}}
{\tiny\left(\!\!
\begin{array}{l}
\text{Askey,}\\
\text{Wilson}\\
\text{1985}
\end{array}\!\!
\right)}
@.\\
@VVV @V{{\text{{\normalsize \bf \em  self-dual orthogonal}}\atop }\atop}
V
{{\text{{\normalsize  \bf\em polynomials}}\atop }\atop}
V@.\\
{\Large
\text{${}_7F_6$}}
{\tiny\left(\!\!
\begin{array}{l}
\text{Dougall}\\
\text{1907}
\end{array}\!\!
\right)}
@>>>
{\Large\text{
${}_8\varphi_7$}}
{\tiny\left(\!\!
\begin{array}{l}
\text{Jackson}\\
\text{1921}
\end{array}\!\!
\right)}
@>>>
{\Large\text{
${}_{10}E_9$}}
{\tiny\left(\!\!
\begin{array}{l}
\text{Frenkel,}\\
\text{Turaev}\\
\text{1997}
\end{array}\!\!
\right)}
\\
@VVV @V{{\text{{\normalsize \bf\em summation formulas}}\atop }\atop }
VV@VVV\\
{\Large\text{
${}_9F_8$}}
{\tiny\left(\!\!
\begin{array}{l}
\text{Wilson}\\
\text{1978}\\
\text{Rahman}\\
\text{1986}
\end{array}\!\!
\right)}
@>>>
{\Large\text{
${}_{10}\varphi_9$}}
{\tiny\left(\!\!
\begin{array}{l}
\text{Rahman}\\
\text{1986}\\
\text{Wilson}\\
\text{1991}
\end{array}\!\!
\right)}
@>>>
\fbox{{\Large${}_{12}E_{11}$}}
{\tiny\left(\!\!
\begin{array}{l}
\text{\cite{sz:cmp}} \\ 
\text{ 2000}\\ 
\text{\cite{spi:theta2}}\\ 
\text{ 2003} \\ 
\end{array}\!\!
\right)}
\\
\end{CD}
$$
\hspace{1.8cm}{\bf \em self-dual biorthogonal rational functions}
\vskip 5mm

\noindent
 expressed in terms of a special $_4\varphi_3$ series
and their argument ``lives" on the most general admissible grid for polynomials.

The next level of complexification of functions indicated in the table refers to the
most general known summation formulas for terminating series of
hypergeometric type. Sequentially, these are the Dougall's $_7F_6$
and Jackson's $_8\varphi_7$ sums going back to the first quarter
of the last century, and the recent result by Frenkel and Turaev
\cite{FT} at the level of $_{10}E_9$ series to be described below.

Finally, until very recent time the most general set of known special functions
satisfying some orthogonality relations and obeying other ``classical''
properties were given by biorthogonal rational functions related to the
very well poised $_9F_8$ and $_{10}\varphi_9$ series. The discrete
measure functions were discovered by Wilson \cite{wil:orthogonal}
and their continuous measure generalizations were derived by Rahman
\cite{rah:integral}. An elliptic extension of the Wilson's
biorthogonal functions with the key self duality property was constructed
by Zhedanov and the author \cite{sz:cmp}. The Rahman's family of
rational functions was lifted to the elliptic level by the author
\cite{spi:theta2}. These functions ``live" on the grids described by the
second order elliptic functions---the most general type of grids for rational
functions admitting a lowering divided difference operator \cite{spi-zhe:pade}.
Moreover, in the elliptic case there appeared even more complicated
objects existing only at this level \cite{spi:theta2}, which go beyond the
space of rational functions of some argument and which satisfy unusual
two index biorthogonality relations.

There exist also non-self-dual three parameter extension of the last row
functions described by the very well poised
${}_9F_8$, ${}_{10}\varphi_9 $, and ${}_{12}E_{11}$ series \cite{sz:cmp},
but many of their properties remain unknown.

In the following we restrict ourselves only to the elliptic hypergeometric
functions and for further details concerning plain and $q$-hypergeometric
objects we refer to the textbooks \cite{aar:special} and \cite{gas-rah:basic},
handbook \cite{koe-swa:askey} and the original papers
\cite{ask-wil:some,rah:integral,wil:orthogonal}. For a description of general
formal unilateral $_sE_r$ and bilateral $_sG_r$ theta hypergeometric series,
see \cite{gas-rah:basic,spi:theta1,spi:thesis}.

\underline{Elliptic hypergeometric integrals}
are described with the help of the bases $q, p$ and
\vskip 1mm
\centerline{
\fbox{$\tilde{q}=e^{-2\pi i \omega_2/\omega_1},\
\tilde{p}=e^{-2\pi i \omega_2/\omega_3},\
r=e^{2\pi i \omega_3/\omega_1},\
\tilde{r}=e^{-2\pi i \omega_1/\omega_3},$} }
\vskip 1mm
\noindent
where $\tilde{q}$, $\tilde{p}$, $\tilde{r}$ are modular transforms
of $q$, $p$, and $r$.

\begin{thm}{\bf (An elliptic analogue of the Meijer function \cite{spi:theta2})}
\par
For incommensurate $\omega_i$ and $|p|,|q|,|r|<1$ general solution of the equations
$$
\Delta(u+\omega_i)=h_i(u)\Delta(u),\quad i=1,2,3,
$$
where
\begin{center}
\begin{tabular}{|c|c|c|c|}
\hline
ell.~fun-s & periods & bases & moduli\\
\hline
\hline
$h_1(u)$ & $\omega_2$, $\omega_3$ & $p$ & $\tau_1=\omega_3/\omega_2$ \\
\hline
$h_2(u)$ & $\omega_1$, $\omega_3$ & $r$ & $\tau_2=\omega_3/\omega_1$ \\
\hline
$h_3(u)$ & $\omega_1$, $\omega_2$ & $q$ & $\tau_3=\omega_1/\omega_2$ \\
\hline
\end{tabular}
\quad $\tau_1=\tau_2\tau_3,$
\end{center}
is:
$$
\Delta(u) =\prod_{j=0}^m
\frac{\mit\Gamma(t_j e^{2\pi i u/\omega_2};p,q)}
{\mit\Gamma(w_j e^{2\pi i u/\omega_2};p,q)} \prod_{j=0}^{m'}
\frac{{\mit\Gamma}(t'_j e^{-2\pi i u/\omega_1};r,\tilde{q})}
{{\mit\Gamma}(w'_j e^{-2\pi i u/\omega_1};r,\tilde{q})}
\ e^{\gamma u}\times \mbox{{\rm constant}},
$$
where
$$
\prod_{j=0}^m t_j=\prod_{j=0}^m w_j,\qquad
\prod_{j=0}^{m'} t'_j=\prod_{j=0}^{m'} w'_j,
$$
and
$$
{\mit\Gamma}(z;p,q)
=\prod_{j,k=0}^\infty
\frac{1-z^{-1}p^{j+1}q^{k+1}}{1-zp^{j}q^{k}},
\quad |q|, |p|<1,
$$
is the standard elliptic gamma function \cite{rui:first}.
\end{thm}

The  function ${\mit\Gamma}(z;p,q)$ satisfies equations
\begin{eqnarray*}
{\mit\Gamma}(z;p,q)&=&{\mit\Gamma}(z;q,p),\quad
{\mit\Gamma}(pq/z;p,q)=1/{\mit\Gamma}(z;p,q),
\\
{\mit\Gamma}(qz;p,q)&=&\theta(z;p){\mit\Gamma}(z;p,q),\\
{\mit\Gamma}(pz;p,q)&=&\theta(z;q){\mit\Gamma}(z;p,q).
\end{eqnarray*}
If we denote $f(u)=\mit\Gamma(e^{2\pi i u/\omega_2};p,q)$, then this function
solves uniquely (up to a multiplicative factor independent on $u$)
the following system of three linear first order finite difference equations
\[\left\{\begin{array}{l}\mbox{\framebox{$f(u+\omega_1)
= \theta(e^{2\pi i u/\omega_2};p)f(u),$}}
\\ f(u+\omega_2) = f(u),\\ f(u+\omega_3)
= \theta(e^{2\pi i u/\omega_2};q)f(u). \end{array}\right.\]
There are two choices of parameters with additional nice properties:
\begin{itemize}
\item [1)]$\gamma=0$ and no $t'_k, w'_k$
(the ``standard" case $|p|,|q|<1$);
\item [2)]$\gamma=0$ and $m'=m$, $t'_j=rt_j$, $w'_j=rw_j$
(the ``unit circle" case).
\end{itemize}

In the second case, gamma function factors combine into the modified
elliptic gamma function introduced in \cite{spi:theta2}:
$$
G(u;\omega_1,\omega_2,\omega_3)
={\mit\Gamma}(e^{2\pi i u/\omega_2};p,q)
{\mit\Gamma}(re^{-2\pi i u/\omega_1};r,\tilde{q}).
$$
This function solves uniquely another system of three equations:
\[\left\{\begin{array}{l}
\mbox{\framebox{$f(u+\omega_1)
= \theta(e^{2\pi i u/\omega_2};p)f(u),$}}
\\
f(u+\omega_2) =\theta(e^{2\pi i u/\omega_1};r) f(u),
\\
  f(u+\omega_3) =e^{-\pi iB_{2,2}(u)} f(u),
\end{array}\right.\]
where
$$
B_{2,2}(u)=\frac{u^2}{\omega_1\omega_2}
-\frac{u}{\omega_1}-\frac{u}{\omega_2}+
\frac{\omega_1}{6\omega_2}+\frac{\omega_2}{6\omega_1}+\frac{1}{2}.
$$
These equations allow us to prove the representation \cite{die-spi:unit}
$$
 G(u;{\bf \omega})=e^{-\pi i P(u)}\Gamma(e^{-2\pi
i \frac{u}{\omega_3}}; \tilde r, \tilde p),
$$
$$
P\left(u+\sum_{n=1}^3\frac{\omega_n}{2}\right)
=\frac{u(u^2-\frac{1}{4}\sum_{k=1}^3\omega_k^2)}{3\omega_1\omega_2\omega_3},
$$
related to modular transformations for the standard
elliptic gamma function \cite{fel-var:elliptic}.
{}From this representation it is easy to see that
$G(u;{\bf \omega})$ is well defined for $|p|,|r|<1$ and $|q|\le1$
(i.e., the $|q|=1$ case is permitted in sharp difference from the
$\Gamma(z;p,q)$ function!).

Permutations $\tilde r\leftrightarrow\tilde p$
and $\omega_1\leftrightarrow\omega_2$ are equivalent. Therefore, we have
$$
G(u;\omega_1,\omega_2,\omega_3)=G(u;\omega_2,\omega_1,\omega_3).
$$
Due to the property $P(\sum_{k=1}^3\omega_k-u)=-P(u)$, we have the
reflection equation
$$
G(a;{\bf \omega})G(b;{\bf \omega})=1,\quad a+b=\sum_{k=1}^3\omega_k.
$$

In the limit $\omega_3\to\infty$, taken in such a way that simultaneously $p,r\to0$,
the modified elliptic gamma function is reduced to the ``unit circle"
$q$-gamma function
$$
\lim_{p,r\to 0} \frac{1}{G(u;{\bf \omega})}
=S(u;\omega_1,\omega_2) = \frac{(e^{2\pi i u/\omega_2}; q)_\infty}
{(e^{2\pi iu/\omega_1}\tilde q; \tilde q)_\infty},
$$
which remains well defined in the limit $|q|\to 1$.
This function appeared in the modern time mathematics in the work of Shintani
\cite{shi:kronecker} as a ratio of Barnes' double gamma functions \cite{bar:theory};
in the works of Faddeev and coauthors \cite{fad:discrete,fkv:strongly}
on the modular double of quantum groups and quantum Liouville theory;
in the work of Jimbo and Miwa \cite{jim-miw:quantum} on solutions of a
$q$-difference equation and related correlation functions in statistical mechanics;
in eigenfunctions of the $q$-Toda chain Hamiltonian \cite{kls:unitary}.
In several independent studies it was named as the double sign function \cite{kur:multiple},
or hyperbolic gamma function \cite{rui:first,rui:generalized}, or non-compact
quantum dilogarithm \cite{fkv:strongly}.
For the operator algebra aspects of this nice function, see \cite{volk}.

\section{Elliptic functions versus balanced, well poised and
 very well poised hypergeometric functions}

Some convenient terminology.

\noindent
\underline{Theta functions}: holomorphic functions $f(x)$ such that
$$
f(x+\omega_2)=e^{ax+b}f(x), \quad
f(x+\omega_3)=e^{cx+d}f(x),
$$
for some $a,b,c,d\in\mathbb{C}$. They have a finite number of zeros in the
parallelogram of periods $\omega_2, \omega_3$, Im$(\omega_2/\omega_3)\neq 0$.
It is not difficult to deduce that
$$
f(x)=e^{P_2(x)}\prod_{j=0}^r \theta_1(x+u_j),
\quad u_j\in {\mathbb C},
$$
for some polynomial of the second order $P_2(x)$.

\noindent
\underline{Meromorphic theta functions}:
ratios of theta functions with different parameters $r, u_j$ and $P_2(x)$.

\noindent
\underline{Elliptic functions}: balanced meromorphic theta functions
$$
f(x)=\prod_{j=0}^r\frac{\theta_1(x+u_j)}{\theta_1(x+v_j)}=
\prod_{j=0}^r \frac{\theta(t_jz;p)}{\theta(w_jz;p)},
$$
where $p=e^{2\pi i \tau},$ $z=e^{2\pi i x},$ $t_j=e^{2\pi i u_j}$,
$w_j=e^{2\pi i v_j}$ with the balancing constraint \fbox{$\prod_{j=0}^r t_j=\prod_{j=0}^r w_j$,}
or \fbox{$\sum_{j=0}^r u_j=\sum_{j=0}^r v_j\, (\mbox{mod}\, 1)$}
guaranteeing that $f(x+1)=f(x)$ and $f(x+\tau)=f(x)$. We can multiply
these functions by arbitrary independent variable $y$ which is omitted
for brevity.

\noindent
\underline{Modular invariant elliptic functions}:
elliptic functions invariant under the action of full
$PSL(2;{\mathbb Z})$ group generated by the relations
$$
f(x;\tau+1)=f(x;\tau),\qquad f(x/\tau; -1/\tau)=f(x;\tau).
$$
Due to the symmetry properties
\begin{eqnarray*}
&& \theta_1(u|\tau+1) = e^{\pi i/4}\theta_1(u|\tau),\qquad
\\
&& \theta_1\Big(\frac{u}{\tau}\Big|-\frac{1}{\tau}\Big)
=-i(-i\tau)^{1/2}e^{\pi iu^2/\tau} \theta_1(u|\tau),
\end{eqnarray*}
elliptic functions are modular if
\fbox{$\sum_{j=0}^r u_j^2 =\sum_{j=0}^r v_j^2\, (\mbox{mod}\, 2\tau).$}
A useful form of the second transformation is
$$
\frac{\theta(e^{-2\pi i\frac{u}{\omega_3}};e^{-2\pi i\frac{\omega_2}{\omega_3}})}
{\theta(e^{2\pi i\frac{u}{\omega_2}};e^{2\pi i\frac{\omega_3}{\omega_2}})}
=i e^{\pi i\frac{\omega_2+\omega_3}{6\omega_2\omega_3}}
e^{\pi i\frac{u^2-u(\omega_2+\omega_3)}{\omega_2\omega_3}}
=e^{\pi iB_{2,2}(u)},
$$
which indicates that the true modular transformation corresponds to the
change $(\omega_2,\omega_3)\to(-\omega_3,\omega_2)$.

\noindent
\underline{Totally elliptic functions:}
elliptic $f(x)$ which are elliptic also in $u_j$, $v_j$
with the same periods. These are elliptic functions with the constraints
\fbox{$v_j=-u_j (\text{mod}\, 1)$} or \fbox{$w_j=t_j^{-1}$} known
in the theory of $q$-hypergeometric series as the
well poisedness conditions. The balancing condition for such
{\em well poised elliptic functions}
$\prod_{j=0}^r t_j=\prod_{j=0}^r t_j^{-1}$ is reduced to
\fbox{$ \prod_{j=0}^r t_j=\pm 1$,} i.e. we have a sign ambiguity!
Totally elliptic functions are invariant under the shifts
$t_j\to p t_j$ $(j=0,1,\ldots,r-1)$ and $z\to pz$. Moreover, they
are automatically modular invariant and satisfy the relation $f(-x)=1/f(x)$
(this relation reduces $y$, the arbitrary multiplier of $f(x)$, to $y=\pm 1$).

We scale now $z\to t_0 z$ and replace parameters $t_j t_0$ by $t_j$
(in particular, we change $t_0^2\to t_0$). As a result, we obtain$$
f(z,\underline{t})=\prod_{j=0}^r \frac{\theta(t_jz;p)}{\theta(t_j^{-1}z;p)}
\to  \mbox{\fbox{ $\displaystyle
\prod_{j=0}^r \frac{\theta(t_jz;p)}{\theta(t_0t_j^{-1}z;p)} $} }
\equiv h(z,\underline{t}).
$$
The balancing condition takes now the form
$\prod_{j=1}^rt_j=\pm t_0^{(r-1)/2}$. Let us take $r=2k+1$ odd and resolve
the sign ambiguity in favor of the relation
\fbox{ $\prod_{j=1}^{2k+1} t_j=+t_0^k$.}
Only for this case there are non-trivial summation and transformation
formulas for series of hypergeometric type. In this case
$h(z,\underline{t})$ is invariant under the shift $t_0\to pt_0$
(accompanied by the compensating transformation $t_{2k+1}\to p^kt_{2k+1}$),
i.e. it is an elliptic function of $\log t_0$
with the same periods as for the $\log z$ variable.
Equivalently, we have
$$
f(p^{1/2}z,p^{1/2}t_0,\ldots,p^{1/2}t_{r-1},p^{-r/2}t_r)
=f(z,t_0,\ldots,t_r),
$$
i.e. there appears interesting symmetry playing with the half
period shifts. We conclude that the total ellipticity
requirement (in appropriate parametrization) fixes the correct form
of the balancing condition in the most interesting
case of odd $r$.

Another important structural constraint leading to interesting
elliptic functions is called the {\em very well poisedness} condition. It
consists in imposing on the well poised elliptic functions
of the restrictions
\par\vskip 1mm
\centerline{\fbox{$
t_{r-3}=q\sqrt{t_0},\ t_{r-2}=-q\sqrt{t_0},\
t_{r-1}=q\sqrt{t_0/p},\ t_{r}=-q\sqrt{pt_0}
$}}
\vskip 1mm\noindent
related to the doubling of the $\theta_1(x)$ function argument.

We call elliptic hypergeometric series and integrals modular, well poised, or
very well poised, if the ratios of their kernels
$c_{n+1}/c_n$ and $\Delta(u+\omega_1)/\Delta(u)$ are modular, well poised,
or very well poised elliptic functions. It is
convenient to introduce special notation for
the very well poised elliptic hypergeometric series \cite{spi:bailey1}:
\begin{eqnarray*}
\lefteqn{
{}_{r+1}E_{r}\bigg({t_0,t_1,\ldots,t_{r-4},q\sqrt{t_0},-q\sqrt{t_0},
q\sqrt{t_0/p},-q\sqrt{pt_0}  \atop
qt_0/t_1,\ldots  ,qt_0/t_{r-4},\sqrt{t_0}, -\sqrt{t_0},\sqrt{pt_0},
-\sqrt{t_0/p}};q,p;-y\bigg)
}\\
&=&\sum_{n=0}^\infty
\frac{\theta(t_0q^{2n};p)}{\theta(t_0;p)}\prod_{m=0}^{r-4}
\frac{(t_m)_n}{(qt_0t_m^{-1})_n}(qy)^n \equiv
{}_{r+1}V_{r}(t_0;t_1,\ldots  ,t_{r-4};q,p;y),
\end{eqnarray*}
where $\prod_{k=1}^{r-4}t_k=\pm t_0^{(r-5)/2}q^{(r-7)/2}$
(for odd $r$ we assume the positive sign, due to the property
described above).
All known applications of these series use a special value
of the argument $y$, $y=1$. Therefore, we shall drop $y$
in the notation of $_{r+1}V_r$ series for this special case.
For $p\to0$, these series reduce to the very well poised
$_{r-1}\varphi_{r-2}$ series denoted by the symbol
$_{r-1}W_{r-2}$ in the monograph \cite{gas-rah:basic}.
Remarkably, the elliptic balancing condition coincides in this case
with the usual balancing condition accepted for these particular
basic hypergeometric series \cite{gas-rah:basic,spi:theta1}.

Various forms of the ellipticity
requirement provide thus an explanation of the origin of the notions
of balancing and very well poisedness for series of hypergeometric
type \cite{spi:theta1}. It is the clarification of these points
that forced the author to change previous notation for
elliptic hypergeometric series \cite{spi:theta1,spi:bailey1}.
In  particular, in this system of conventions accepted in
\cite{gas-rah:basic,ros:elementary,spi:theta2}, etc the symbol
$_{r+1}E_r$ used in the papers \cite{die-spi:elliptic,kmnoy,sz:cmp}
should read as $_{r+3}E_{r+2}$ or $_{r+3}V_{r+2}$.

If we take $r=9$, $t_4=q^{-N} (N\in\N),$ $\prod_{m=1}^5 t_m=qt_0^2$,
$y=1$, then
\vskip 1mm \centerline{\framebox{$\displaystyle
{}_{10}V_{9}(t_0;t_1,\ldots  ,t_{5};q,p)=\frac{(qt_0)_N(\frac{qt_0}{t_1t_2})_N
(\frac{qt_0}{t_1t_3})_N(\frac{qt_0}{t_2t_3})_N }
{(\frac{qt_0}{t_1t_2t_3})_N(\frac{qt_0}{t_1})_N
(\frac{qt_0}{t_2})_N(\frac{qt_0}{t_3})_N}.
$}} \vskip 1mm \noindent
This is the Frenkel-Turaev summation formula \cite{FT} (for
its elementary proofs, see, e.g., \cite{ros:elementary,spi-zhe:theory}),
which is reduced in the limit $p\to 0$ to the Jackson sum for terminating
very well poised balanced ${}_8\varphi_7$ series.

\section{The univariate elliptic beta integral}

The elliptic beta integral is
the simplest very well poised elliptic hypergeometric integral.

\begin{thm}{\bf (The standard elliptic beta integral \cite{spi:beta}) }
\par
Let $t_1,\ldots,t_6\in \C$, $|t_j|<1$, $\prod_{j=1}^6 t_j=pq$, and $|p|,|q|<1$.
Then
\vskip 1mm \centerline{\framebox{$\displaystyle
\kappa\int_{{\mathbb T}}
\frac
{\prod_{k=1}^6{\mit\Gamma}(t_kz;p,q){\mit\Gamma}(t_kz^{-1};p,q)}
{{\mit\Gamma}(z^2;p,q){\mit\Gamma}(z^{-2};p,q)}
\frac{dz}{z}
=\prod_{1\le j<k\le 6}{\mit\Gamma}(t_jt_k;p,q),
$}} \vskip 1mm \noindent
where ${\mathbb T}$ is the positively oriented unit circle $|z|=1$ and
\vskip 1mm \centerline{\framebox{$\displaystyle
\kappa=\frac{(q;q)_\infty(p;p)_\infty}{4\pi i}.
$}} \vskip 1mm
\end{thm}

The first proof of this integration formula used an elliptic generalization of the
Askey's method \cite{ask:beta} which required some contiguous relations for
the left-hand side expression and Bailey's $_2\psi_2$ summation formula.
A very simple proof has been found later on in \cite{spi:short}.

The elliptic beta integral is the most general univariate
beta type integral found so far. It serves as a measure in the biorthogonality
relations for a particular system of functions to be described below.
After taking the limit  $p\to 0$, our integral is reduced to the Rahman's
$q$-beta integral \cite{rah:integral}
\begin{eqnarray*}
&& \frac{(q;q)_\infty}{4\pi i}\int_{{\mathbb T}}
\frac{(z^2;q)_\infty(z^{-2};q)_\infty(Az;q)_\infty(Az^{-1};q)_\infty}
{\prod_{m=1}^5(t_mz;q)_\infty(t_mz^{-1};q)_\infty}
\frac{dz}{z}
\\ && \makebox[4em]{}
=\frac{\prod_{m=1}^5(At_m^{-1};q)_\infty}{\prod_{1\le j<k\le 5}
(t_jt_k;q)_\infty},
\end{eqnarray*}
where $A=\prod_{m=1}^5t_m$, $|t_m|<1$. This integral determines the measure
for Rahman's family of continuous biorthogonal rational functions
\cite{rah:integral}.

If we take now the limit $t_5\to 0$, then
we obtain the celebrated Askey-Wilson integral
$$
\frac{(q;q)_\infty}{4\pi i}\int_{{\mathbb T}}
\frac{(z^2;q)_\infty(z^{-2};q)_\infty}
{\prod_{m=1}^4(t_mz;q)_\infty(t_mz^{-1};q)_\infty}
\frac{dz}{z}
=\frac{(t_1t_2t_3t_4;q)_\infty}{\prod_{1\le j<k\le 4}(t_jt_k;q)_\infty},
$$
determining the measure in orthogonality relations for the most general
set of classical orthogonal polynomials \cite{ask-wil:some}.

Careful analysis of the structure of residues of the integrand's poles
allows one to deduce the Frenkel-Turaev summation formula
out of the elliptic beta integral \cite{die-spi:elliptic}.
We suppose that $|t_m|<1,\, m=1,\ldots ,4,$ $|pt_5|<1<|t_5|$, $|pq|<|A|$,
$A=\prod_{s=1}^5t_s$, and
assume also that the arguments of all $t_s,\, s=1,\ldots,5,$
and $p,q$ are linearly independent over $\mathbb{Z}$. We denote
$C$ a contour separating sequences of integrand's poles
at $z=t_sq^jp^k$ and $A^{-1}q^{j+1}p^{k+1}$,
from their reciprocals at $z=t_s^{-1}q^{-j}p^{-k},$
$Aq^{-j-1}p^{-k-1}$, $j,k\in \mathbb{N}$.
Then we obtain the following residue formula:
$$
\kappa\int_C \Delta_E(z,\underline{t})\frac{d z}{z} =
\kappa\int_\mathbb{T} \Delta_E(z,\underline{t})\frac{d z}{z}
+c_0(\underline{t}) \sum_{\stackrel{n \geq 0}{|t_0q^n|>1}}
\nu_n(\underline{t}),
$$
with
$$
\Delta_E(z,\underline{t})
=\frac{\prod_{m=1}^5\Gamma(t_mz^\pm;p,q)}{\Gamma(z^{\pm2};p,q)\Gamma(Az^\pm;p,q)},
$$
$\Gamma(az^\pm;p,q)\equiv\Gamma(az;p,q)\Gamma(az^{-1};p,q)$, and
\begin{eqnarray*} \nonumber
&&c_0(\underline{t}) =
\frac{\prod_{m=1}^4\Gamma(t_mt_5^\pm;p,q)}
{\Gamma(t_5^{-2};p,q)\Gamma(A t_5^\pm;p,q)}, \\
&&\nu_n(\underline{t}) = q^n\,
 \frac{\theta(t_5^2q^{2n};p)}{\theta(t_5^2;p)}
\prod_{m=0}^5 \frac{(t_mt_5)_n}
     {(qt_m^{-1}t_5)_n},
\end{eqnarray*}
where we have introduced a new parameter $t_0$ via the relation
$\prod_{m=0}^5t_m $ $=q$.
In the limit $t_5t_4\to q^{-N},\, N\in\mathbb{N}$,
values of the integral on the left-hand side of this formula and of the factor
$c_0(\underline{t})$ in front of the residues sum on the right-hand side
blow up, but the integral over the unit circle $\mathbb{T}$ remains finite.
Dividing all the terms by $c_0(\underline{t})$ and taking the limit, we obtain
the summation formula presented in the end of the previous section.

Using the modified elliptic gamma function it is not difficult to deduce
out of the standard elliptic beta integral its ``unit circle'' analogue
remaining well defined for $|q|=1$.

\begin{thm}{\bf (The modified  elliptic beta integral \cite{die-spi:unit}) }
\par
We suppose that Im$(\omega_1/\omega_2)\geq 0$ and Im$(\omega_3/\omega_1)>0$,
Im$(\omega_3/\omega_2)>0$ and $g_j\in\C$, $j=1,\ldots,6$,
Im$(g_j/\omega_3)<0,$ together with the constraint
$\sum_{j=1}^6g_j=\sum_{k=1}^3\omega_k.$ Then
$$
\tilde\kappa\int_{-\omega_3/2}^{\omega_3/2} \frac{\prod_{j=1}^6 G(g_j\pm
u;{\bf \omega})} {G(\pm 2u;{\bf \omega})} \frac{du}{\omega_2}
= \prod_{1\leq j<m\leq 6}G(g_j+g_m;{\bf \omega}),
$$
where
$$
\tilde\kappa= -\frac
{(q;q)_\infty(p;p)_\infty(r;r)_\infty}{2(\tilde q;\tilde q)_\infty}.
$$
Here the integration is taken along the cut
with the end points $-\omega_3/2$ and $\omega_3/2$. We use also the
convention that $G(a\pm b;{\bf \omega})\equiv G(a+b;{\bf \omega})G(a-b;
{\bf \omega})$.
\end{thm}

If we take Im$(\omega_3)\to\infty$ in such a way that $p,r\to 0$,
then this integral reduces to a Mellin-Barnes type $q$-beta integral.
More precisely, for $\omega_{1,2}$ such that Im$(\omega_1/\omega_2)\geq 0$
and Re$(\omega_1/\omega_2)> 0$, we substitute $g_6=\sum_{k=1}^3\omega_k
-\mathcal{A}$, where $\mathcal{A}=\sum_{j=1}^5g_j$ and apply the inversion
formula for $G(u;{\bf \omega})$. Then we set $\omega_3=it\omega_2$,
$t\to +\infty$, and obtain formally
$$
\int_{\mathbb{L}}\frac{S(\pm 2u, \mathcal{A}\pm
u;{\bf\omega})} {\prod_{j=1}^5 S(g_j\pm
u;{\bf \omega})}\frac{du}{\omega_2}= -2\frac{(\tilde
q;\tilde q)_\infty}{(q;q)_\infty}
\frac{\prod_{j=1}^5S(\mathcal{A}-g_j;{\bf \omega})} {
\prod_{1\leq j<m\leq 5}S(g_j+g_m;{\bf \omega})},
$$
where the integration is taken along the line
$\mathbb{L}\equiv i\omega_2\mathbb{R}$.
Here parameters are subject to the constraints Re$(g_j/\omega_2)>0$ and
Re$((\mathcal{A}-\omega_1)/\omega_2)<1$.
This integral was rigorously proven first in \cite{sto:hyperbolic}
and a quite simple proof was given in \cite{spi:short} in a more general
setting.

\section{An elliptic analogue of the ${}_2F_1$ function}

We consider the double integral
$$\kappa \int_{C^2}
\frac
{\prod_{j=1}^3{\mit\Gamma}(a_jz^\pm, b_jw^\pm)\;{\mit\Gamma}(cz^\pm w^\pm)}
{{\mit\Gamma}(z^{\pm 2},w^{\pm 2},c^2Az^{\pm},c^2Bw^{\pm})}
\frac{dz}{z}\frac{dw}{w},
$$
where $a_j,b_j,c\in {\mathbb C}$, $A=a_1a_2a_3$, $B=b_1b_2b_3$, and
$C$ is a contour separating converging to zero sequences
of poles in $z$ and $w$ from the diverging ones, and
$$
{\mit\Gamma}(t_1,\ldots  ,t_k)\equiv
{\mit\Gamma}(t_1;p,q)\cdots{\mit\Gamma}(t_k;p,q).
$$
Applying the elliptic beta integral formula to integrations with respect to $z$
or $w$ (permutation of integrations is allowed since the integrand is bounded),
we obtain a symmetry transformation for a pair of
elliptic hypergeometric integrals \cite{spi:theta2}
\begin{eqnarray*}
&& \prod_{j=1}^3
\frac{{\mit\Gamma}(A/a_j)}{{\mit\Gamma}(c^2A/a_j)}
\int_{C} \frac
{\prod_{j=1}^3{\mit\Gamma}(ca_jz^\pm, b_jw^\pm)}
{{\mit\Gamma}(z^{\pm 2},cAz^{\pm},c^2Bw^{\pm})}
\frac{dz}{z}\\
&& \makebox[2em]{}
= \prod_{j=1}^3
\frac{{\mit\Gamma}(B/b_j)}{{\mit\Gamma}(c^2B/b_j)}
\int_{C} \frac {\prod_{j=1}^3{\mit\Gamma}(a_jz^\pm, cb_jw^\pm)}
{{\mit\Gamma}(z^{\pm 2},c^2Az^{\pm},cBw^{\pm})}
\frac{dz}{z}.
\end{eqnarray*}
This is an elliptic analogue of the four term Bailey transformation for
non-terminating ${}_{10}\varphi_{9}$ series.
It cannot be written yet as some relation for infinite $_{12}V_{11}$
elliptic hypergeometric series due to the severe problems with their
convergence at the boundary values of the argument $|y|=1$.

We denote $t_{1,2,3}=ca_{1,2,3},$ $t_4=pq/cA,$  $t_{5,6,7}=b_{1,2,3},$
$t_8=pq/c^2B$ and introduce the elliptic hypergeometric function---an
elliptic analogue of the Gauss hypergeometric function
\vskip 1mm
\centerline{\framebox{
$\displaystyle
V(\underline{t};p,q)=\kappa\int_C\frac{\prod_{j=1}^8
{\mit\Gamma}(t_jz^{\pm})}{{\mit\Gamma}(z^{\pm 2})}\frac{dz}{z},
\qquad  \prod_{j=1}^8 t_j=p^2q^2.
$}}
\par\vskip 2mm \noindent
Due to the reflection equation for $\Gamma(z;p,q)$ function, we have
$$
V(\underline{t};p,q)\Big|_{t_7t_8=pq}=\prod_{1\leq j<k\leq6}\Gamma(t_jt_k;p,q),
$$
which is the elliptic beta integration formula (evidently, in this
relation $t_7$ and $t_8$ can be replaced by any other pair of parameters).

In the notation $V(\underline{t})=V(\underline{t};p,q)$, the transformation
derived above reads
\vskip 1mm
(i) \framebox{
$\displaystyle
V(\underline{t})=\prod_{1\le j<k\le 4}{\mit\Gamma}(t_jt_k,t_{j+4}t_{k+4})\,
V(\underline{s}),
$} \par\vskip 1mm \noindent
where
\[
\left\{
\begin{array}{cl}
s_j =\varepsilon^{-1} t_j&   ( j=1,2,3,4)  \\
s_j = \varepsilon t_j &    (j=5,6,7,8)
\end{array}
\right.
\quad \varepsilon=\sqrt{\frac{t_1t_2t_3t_4}{pq}}=\sqrt{\frac{pq}{t_5t_6t_7t_8}}.
\]
\par \vspace{1mm}

We repeat this transformation with $s_3,s_4,s_5,s_6$ playing the role of
$t_1,t_2,t_3,t_4$ and permute parameters $t_3,t_4$ with $t_5,t_6$
in the result. This yields

(ii) \framebox{
$\displaystyle
V(\underline{t})=\prod_{j,k=1}^4
{\mit\Gamma}(t_jt_{k+4})\ V(T^{1\over 2}\!/t_1,\ldots,T^{1\over 2}\!/t_4,
U^{1\over 2}\!/t_5,\ldots,U^{1\over 2}\!/t_8),
$}
\par\vskip 1mm\noindent
where
$ T=t_1t_2t_3t_4$ and $ U=t_5t_6t_7t_8. $
\par \vspace{1mm}
We equate now the right-hand sides of relations  (i) and (ii), express $t_j$ parameters
in terms of $s_j$ and obtain
\vskip 1mm
(iii) \framebox{
$\displaystyle
V(\underline{s})=\prod_{1\le j<k\le 8}{\mit\Gamma}(s_js_k)\,
V(\sqrt{pq}/\underline{s}),
$}  \par\vskip 1mm\noindent
where $\sqrt{pq}/\underline{s}=(\sqrt{pq}/s_1,\ldots,\sqrt{pq}/s_8)$.

Transformations (ii) and (iii) were proven by Rains \cite{rai:trans}
in a strai\-ghtforward manner using evaluations of determinants of theta functions
on a dense set of parameters. However, as we just have seen \cite{spi:thesis}, they
are mere repetitions of the key transformation (i).

It is convenient to set temporarily $t_j=e^{2\pi i x_j}(pq)^{1/4}$. We
take vectors $x\in\R^8$ and denote as $x=\sum_{i=1}^8x_ie_i$ their
standard decomposition in the orthonormal basis $e_i$,
$\langle e_i,e_j\rangle=\delta_{ij}$. Then the balancing condition
implies $\sum_{i=1}^8 x_i=0$ which defines a hyperplane $Y$ orthogonal to
the vector $e_1+\ldots+e_8$. Considering reflections
$x\to x-{2\langle v, x\rangle}\, v/{ \langle v, v\rangle}$
with respect to the hyperplane normal to some vector $v\in Y$, it is not difficult to see
that the transformation of coordinates  in (i) corresponds to the
reflection  with respect to the vector
$v=(\sum_{i=5}^8e_i-\sum_{i=1}^4 e_i)/2$,
which has the canonical normalization of the length $\langle v, v\rangle=2$.

The elliptic hypergeometric function $V(\underline{t})$ appeared for the
first time in our paper \cite{spi:theta2} together with the transformation (i).
However, it was not recognized there that (i) and permutations of parameters
$t_i\leftrightarrow t_j$  generate the exceptional $E_7$ Weyl group of
symmetries: the function
$V(\underline{t})/\prod_{1\leq k<l\leq 8} \sqrt{\Gamma(t_kt_l)}$ is simply
invariant under these transformations. This fact was understood
at the level of series in \cite{kmnoy} (where, actually, only
the $E_6$ group is valid since one of the parameters is fixed to
terminate the series) and for general function $V(\underline{t})$
in \cite{rai:trans}.

For elliptic hypergeometric functions it is convenient to keep two systems
of notation---the ``multiplicative" system, described above, and the ``additive"
one \cite{gas-rah:basic,spi:theta1,spi:bailey1}. Therefore
 we define the function
$$
v(\underline{g};\omega_1,\omega_2,\omega_3)\equiv
V(e^{2\pi i g_1/\omega_2},\ldots,e^{2\pi i g_8/\omega_2};
e^{2\pi i\omega_1/\omega_2},e^{2\pi i\omega_3/\omega_2}),
$$
where  $\sum_{j=1}^8g_k=2\sum_{k=1}^3\omega_k$.
It will be useful for a description of elliptic
hypergeometric equation solutions.

\section{Contiguous relations and the elliptic hypergeometric equation}

The fundamental addition formula for elliptic theta functions can be
written in the following form
\par \vskip 1mm \fbox{$\displaystyle
\theta\left(xw,\frac{x}{w},yz,\frac{y}{z};p\right)
-\theta\left(xz,\frac{x}{z},yw,\frac{y}{w};p\right)
= \frac{y}{w}\theta\left(xy,\frac{x}{y},wz,\frac{w}{z};p\right),
$} \par\vskip 2mm \noindent
where $w,x,y,z$ are arbitrary complex variables.
If we denote $y=t_6, w=t_7, $ and $x=q^{-1}t_8$, then this identity for
theta functions is equivalent to the following $q$-difference equation
\begin{eqnarray*}
&& \Delta(z,t_1,\ldots,t_5,qt_6,t_7,q^{-1}t_8)
-\frac{\theta(t_6t_7^\pm;p)}{\theta(q^{-1}t_8t_7^\pm;p)} \Delta(z,\underline{t})
\\ && \makebox[2em]{}
=\frac{t_6}{t_7}\frac{\theta(q^{-1}t_8t_6^\pm;p)}{\theta(q^{-1}t_8t_7^\pm;p)}
\Delta(z, t_1,\ldots,t_6, qt_7,q^{-1}t_8),
\end{eqnarray*}
where $\Delta(z,\underline{t})=\prod_{k=1}^8\Gamma(t_kz^\pm)/\Gamma(z^{\pm2})$
is the $V$-function integrand.
Integrating now this equality over $z$ along the contour $C$, we
derive the first contiguous relation
\begin{eqnarray*} &&
t_7\theta\left(t_8t_7/q,t_8/qt_7;p\right)V(qt_6,q^{-1}t_8)
-(t_6\leftrightarrow t_7)
 \\ && \makebox[4em]{}
=t_7\theta\left(t_6t_7,t_6/t_7;p\right) V(\underline{t}),
\end{eqnarray*}
which was used in the first proof of the elliptic beta integral
\cite{spi:beta}. Here $V(qt_6,q^{-1}t_8)$ denotes $V(\underline{t})$ with
the parameters $t_6$ and $t_8$ replaced by $qt_6$ and $q^{-1}t_8$ respectively and
$(t_6\leftrightarrow t_7)$ means permutation of the parameters in the
preceding expression.

In the same way as in the case of series \cite{spi-zhe:theory},
we can substitute symmetry transformation (iii) of the previous section
into this equation and obtain the second contiguous relation
\begin{eqnarray*}
&& t_6\theta(t_7/qt_8;p)\prod_{k=1}^5\theta(t_6t_k/q;p)V(q^{-1}t_6,qt_8)
-(t_6\leftrightarrow t_7)\\
&&  \makebox[4em]{}
=t_6\theta(t_7/t_6;p)\prod_{k=1}^5\theta(t_8t_k;p) V(\underline{t}).
\end{eqnarray*}
An appropriate combination of these two equations yields
$$
b(\underline{t})\Big(U(qt_6,q^{-1}t_7)-U(\underline{t})\Big)
+(t_6\leftrightarrow t_7)+ U(\underline{t})=0,
$$
where
$$
U(\underline{t})=\frac{V(\underline{t})}
{\prod_{k=1}^7\mit\Gamma(t_kt_8,t_k/t_8)}
$$
and the potential
\begin{eqnarray*}
&& b(\underline{t})=\frac{\theta(t_6/qt_8,t_6t_8,t_8/t_6;p)}
                 {\theta(t_6/t_7,t_7/qt_6,t_6t_7/q;p)}
\prod_{k=1}^5\frac{\theta(t_7t_k/q;p)}{\theta(t_8t_k;p)}
\\ && \makebox[2em]{}
= \frac{\theta(qt_0/t_6,t_0t_6,t_0/t_6;p)}{\theta(t_6/t_7,qt_6/t_7,q/t_6t_7;p)}
\prod_{k=1}^5\frac{\theta(q/t_7t_k;p)}{\theta(t_0t_k;p)}
\end{eqnarray*}
(the second expression is obtained after setting $t_8=p^2t_0$)
is a modular invariant elliptic function of variables $g_1,\ldots,g_7$
($t_j=e^{2\pi i g_j/\omega_2}$).

If we substitute $t_6=az, t_7=a/z$ and replace $U(\underline{t})$ by
some unknown function $f(z)$, then  we obtain a $q$-difference equation of the
second order called the {\em elliptic hypergeometric equation}:
\begin{eqnarray*}
&& \frac{\theta(az/qt_8,at_8z,t_8/az;p)}
{\theta(z^2,1/qz^2;p)}\prod_{k=1}^5
\theta(at_k/qz;p)\left( f(qz)-f(z)\right)
\\ \nonumber && \makebox[2em]{}
+ \frac{\theta(a/qt_8z,at_8/z,t_8z/a;p)}
{\theta(1/z^2,z^2/q;p)}\prod_{k=1}^5
\theta(at_kz/q;p)\left( f(q^{-1}z)-f(z)\right)
\\  && \makebox[8em]{}
+\theta(a^2/q;p) \prod_{k=1}^5\theta(t_kt_8;p)\, f(z)=0,
\end{eqnarray*}
where $t_8=p^2q^2/a^2\prod_{k=1}^5t_k$. We have found already one functional
solution of this equation $U(\underline{t})$ in the restricted region of parameters.
The second independent solution
can be obtained after scaling any of the parameters $a,t_1,\ldots,t_5$ or $z$ by $p$.
We can replace also the
standard elliptic gamma functions in the definition of $U(\underline{t})$
by the modified elliptic gamma functions and get new solutions of the
elliptic hypergeometric equation.
Indeed, we can rewrite the elliptic hypergeometric equation in the
``additive" notation $t_j=e^{2\pi ig_j/\omega_2}$.
Then the function
$$
v^{mod}(\underline{g};{\bf \omega})
= \int_{-\omega_3/2}^{\omega_3/2} \frac{\prod_{j=1}^8
G(g_j\pm x;{\bf\omega})}{G(\pm 2x;{\bf\omega})}\frac{dx}{\omega_2},
$$
where $\sum_{j=1}^8g_j=2\sum_{k=1}^3\omega_k$,
defines its solution linearly independent from $V(\underline{t})$,
provided we impose appropriate restrictions upon the parameters. Namely, we
should line
up sequences of the integrand's poles to the left or right of the line
passing through the points $-\omega_3/2$ and $\omega_3/2$.
Evidently, $E_7$ symmetry remains intact which follows from the fact
that in the derivation of relevant properties of the $V(\underline{t})$ function
we used only the first (boxed) equation for the elliptic gamma function
$\Gamma(z;q,p)$ which coincides with one of the equations for $G(u;{\bf \omega})$.
Simple computations yield the relation
$$
v^{mod}(\underline{g};{\bf \omega})=
\frac{2\omega_3e^{2\pi i (P(0)-\sum_{j=1}^8P(g_j))}}
{\omega_2(\tilde p;\tilde p)_\infty(\tilde r;\tilde r)_\infty}\,
v(\underline{g};\omega_1,-\omega_3,\omega_2),
$$
showing that this solution is proportional to the modular transformation
of the function $v(\underline{g};\omega_1,\omega_2,\omega_3)$.

Now we shift $g_{7,8}\to g_{7,8}+\sum_{k=1}^3\omega_k$ and take the limit
Im$(\omega_3)\to\infty$ in such a way that $p,r\to 0$. Then our
$v^{mod}$-function is reduced to
$$
s(\underline{g};\omega_1,\omega_2)
= \int_{\L} \frac{S(\pm 2u, -g_7\pm u,-g_8\pm u;{\bf\omega})}
{\prod_{j=1}^6 S(g_j\pm u;{\bf\omega})}\frac{du}{\omega_2},
$$
where $\sum_{j=1}^8g_j=0$.
This is a $q$-hypergeometric function which is well defined for $|q|=1$
and which provides a functional solution of the $p=0$ degeneration of the
elliptic hypergeometric equation.

It should be noticed that $V(\underline{t})$
satisfies not one, but much more equations of the derived type due to
the permutational symmetry in all its parameters, including the equation obtained
after the permutation of $q$ and $p$. Most probably there is only one
function satisfying all of them, since the linearly independent solutions
break one of its symmetries, $E_7$ or $p\leftrightarrow q$.

At the level of $q$-hypergeometric functions, in the limit $p\to 0$
we obtain the equation investigated in detail by Gupta and Masson
\cite{gup-mas:contiguous}. They derived its functional solutions in the
form of  special combinations of non-terminating $_{10}\varphi_9$
series, the integral representation for which has been found earlier by
Rahman \cite{rah:integral} and to which our representation for $V(\underline{t})$
is reduced in the limit $p\to 0$.

In a similar way one can construct contiguous relations for elliptic
${}_{12}V_{11}$ series with $y=1$. Denoting
${\cal E}(\underline{t}) \equiv {_{12}}V_{11}(t_0;t_1,\ldots,t_7;q,p),$
where $\prod_{m=1}^7t_m=t_0^3q^2$ and $t_m=q^{-n},\, n\in\mathbb{N},$ for
some $m$, we have the first relation \cite{sz:cmp,spi-zhe:theory}
\begin{eqnarray*}
&&
{\cal E}(\underline{t}) - {\cal E}(q^{-1}t_6,
qt_7) =\frac{\theta(qt_0,q^2t_0,qt_7/t_6,t_6t_7/qt_0;p)}
{\theta(qt_0/t_6,q^2t_0/t_6,t_0/t_7,t_7/qt_0;p)}  
\\ && \makebox[2em]{} \times
\prod_{r=1}^5\frac{\theta(t_r;p)}{\theta(qt_0/t_r;p)}\,
{\cal E}(q^2t_0;qt_1,\ldots,qt_5,t_6,qt_7),
\end{eqnarray*}
and the second one
\begin{eqnarray*}
\lefteqn{\frac{\theta(t_7;p)\prod_{r=1}^5 \theta(t_rt_6/qt_0;p)}
{\theta(t_6/qt_0,t_6/q^2t_0,t_6/t_7;p)}\,
{\cal E}(q^2t_0;qt_1,\ldots,qt_5,t_6,qt_7) } &&
\\ && \makebox[2em]{}
+\frac{\theta(t_6;p)\prod_{r=1}^5 \theta(t_rt_7/qt_0;p)}
{\theta(t_7/qt_0,t_7/q^2t_0,t_7/t_6;p)}\,
{\cal E}(q^2t_0;qt_1,\ldots,qt_6,t_7)
\nonumber \\ && \makebox[6em]{}
=\frac{\prod_{r=1}^5\theta(qt_0/t_r;p)}{\theta(qt_0,q^2t_0;p)}
\, {\cal E}(\underline{t}).
\end{eqnarray*}

These relations can also be obtained after application of the residue
calculus  similar to the one described above.
For this it is necessary to take one of the parameters of $V(\underline{t})$
outside of the contour $C$ and represent this elliptic hypergeometric
function as a sum of an integral
over $C$ and of the residues picked up during this procedure. An accurate limit
for one of the parameters converting the sum of residues into the terminating
$_{12}V_{11}$ series brings in the needed contiguous relations, which
take the described form after changing notation.

An appropriate combination of these two relations yields
\begin{eqnarray*}\nonumber
\lefteqn{
\frac{\theta(t_6,t_0/t_6,qt_0/t_6;p)}
{\theta(qt_6/t_7,t_6/t_7;p)}\prod_{r=1}^5\theta(qt_0/t_7t_r;p)
\left({\cal E}(qt_6,q^{-1}t_7)-
{\cal E}(\underline{t})\right)
 } &&
\\  \nonumber
&& + \frac{\theta(t_7,t_0/t_7,qt_0/t_7;p)}{\theta(qt_7/t_6,t_7/t_6;p)}
\prod_{r=1}^5\theta(qt_0/t_6t_r;p)\left(
{\cal E}(q^{-1}t_6,qt_7) -
{\cal E}(\underline{t}) \right)\\
&& \makebox[4em]{}
+\theta(qt_0/t_6t_7;p)\prod_{r=1}^5\theta(t_r;p)\,
{\cal E}(\underline{t})=0,
\end{eqnarray*}
which is another form of the elliptic hypergeometric equation.

\section{Applications in mathematical physics}

The theory outlined above did not emerge from scratch. It appeared from
long time developments
in mathematical physics related to classical and quantum completely
integrable systems.
Below we list some of the known applications of elliptic hypergeometric
series and integrals.

\begin{itemize}

\item[(1)]
Elliptic solutions of the Yang-Baxter equation (elliptic $6j$-sym\-bols)
sequentially derived by Baxter \cite{B}, Andrews, Baxter and Forrester \cite{ABF},
Date, Jimbo, Kuniba, Miwa, and Okado \cite{DJKMO} appear to combine into
terminating ${}_{12}V_{11}$
series with special discrete values of parameters, as it was shown
by Frenkel and Turaev in their profound paper \cite{FT}. For a recent work
in this direction including the algebraic aspects of the elliptic $6j$-symbols,
see \cite{konno,rai:abelian,ros:elementary,ros:sklyanin}.
Since solvable two-dimensional statistical mechanics models are related to
the conformal field theory \cite{bm,zub}, it is natural to expect that elliptic
hypergeometric functions will emerge there as well.

\item[(2)] In a joint work with Zhedanov \cite{sz:cmp}, the
terminating ${}_{12}V_{11}$ series with arbitrary continuous
parameters  were discovered as solutions of the linear problem for some
classical integrable system. More precisely, these series emerged from self-similar
solutions of the discrete time chain associated with biorthogonal
rational functions which generalizes ordinary and
relativistic discrete-time Toda chains.

\item[(3)] As shown by Kajiwara, Masuda, Noumi, Ohta, and Yamada \cite{kmnoy},
Sakai's elliptic Painlev\'e equation
\cite{sak} has a solution expressed in terms of the terminating ${}_{12}V_{11}$
series. This observation follows from the reduction of corresponding nonlinear
second order finite difference equation to the elliptic hypergeometric equation.
Therefore, $V(\underline{t})$ also provides its solution.
Moreover, the function $v(\underline{g};\omega_1,-\omega_3,\omega_2)$,
well defined in the $|q|=1$ region, plays a similar role  \cite{spi:thesis}
since it defines
an independent solution of the elliptic hypergeometric equation with the $E_7$
symmetry. More complicated solutions of this equation expressed
in terms of the multiple elliptic hypergeometric integrals were presented by Rains
at this workshop \cite{rai:talk}.

\item[(4)] Elliptic hypergeometric functions provide particular solutions
of the finite difference (relativistic) analogues of the elliptic
Calogero-Sutherland type models \cite{spi:thesis}.
This application is outline below and in the last section.

\end{itemize}

The original investigations of completely integrable many particles systems
on the line (or circle) by Calogero, Sutherland and Moser were continued
by Olshanetsky and Perelomov \cite{OP} who showed that such models
are naturally associated with the root systems.
Relativistic (or finite-difference) generalizations
of these models have been discovered by Ruijsenaars \cite{rui:complete}
who worked out the $A_n$ root system case in detail.
The corresponding eigenvalue problem is also known to be related to the
Macdonald polynomials \cite{mac:constant}.
Inozemtsev \cite{ino:lax} has investigated the most general $BC_n$ root
system extension of the Heun equation absorbing previously derived
differential operator models. In a further step, van Diejen
\cite{die:integrability} has unified Inozemtsev and Ruijsenaars models by
coming up with even more general integrable model, which was investigated
in detail by Komori and Hikami \cite{kom-hik:quantum}. A special degeneration
of this model to the trigonometric level corresponds to
the Koornwinder polynomials \cite{koo:pol}.

The Hamiltonian of the van Diejen model has the form
$$
{\cal H}=\sum_{j=1}^n\Big(A_j(\underline{z})T_j
+A_j(\underline{z}^{-1})T_j^{-1}\Big)+u(\underline{z}),
$$
where $u(\underline{z})$ is some complicated explicit combination of theta functions,
$T_jf(\ldots,z_j,\ldots)=f(\ldots,qz_j,\ldots ),$ and
$$
 A_j(\underline{z})=\frac{\prod_{m=1}^8\theta(t_mz_j;p)}{\theta(z_j^2,qz_j^2;p)}
\prod_{k=1\atop \,\ \ne j}^n
\frac{\theta(tz_jz_k,tz_jz_k^{-1};p)}{\theta(z_jz_k,z_jz_k^{-1};p)}.
$$
If we impose the constraint $t^{2n-2}\prod_{m=1}^8t_m=p^2q^2$, then
the operator ${\cal H}$ can be rewritten in the form
$$
{\cal D}=\sum_{j=1}^n\Big(A_j(\underline{z})(T_j-1)
+A_j(\underline{z}^{-1})(T_j^{-1}-1)\Big)
$$
up to some additive constant independent on variables $z_j$ (for details, see
\cite{die:integrability,kom-hik:quantum,rui:cadiz}).

The standard eigenvalue problem,
$\mathcal{D}f(\underline{z})=\lambda f(\underline{z}),$
in the univariate case $n=1$ looks like
\begin{eqnarray*}
&& \frac{\prod_{j=1}^8\theta(t_jz;p)}{\theta(z^2,qz^2;p)} (f(qz)-f(z))
\\ && \makebox[2em]{}
+ \frac{\prod_{j=1}^8\theta(t_jz^{-1};p)}{\theta(z^{-2},qz^{-2};p)}
(f(q^{-1}z)-f(z))=\lambda f(z).
\end{eqnarray*}
Comparing it with the elliptic hypergeometric equation in the form derived
in \cite{spi:theta2}, which will be described in the next section, we see
that they coincide for a restricted choice of parameters $t_6=t_5/q$
and a special eigenvalue for the Hamiltonian $\mathcal{D}$,
$\lambda=-\kappa_\mu$ (a similar observation has been done by Komori).

However, connections between the elliptic hypergeometric functions and
Calogero-Sutherland type models are deeper than it is just indicated. Let us introduce
the inner product
$$
\langle \varphi,\psi \rangle=\kappa
\int_C  \frac {\prod_{m=1}^8{\mit\Gamma}(t_mz^{\pm})}
{{\mit\Gamma}(z^{\pm 2})} \,\varphi(z)\,\psi(z)\, \frac{dz}{z},
$$
where contour $C$ separates sequences of the kernel's poles converging to
$z=0$ from those diverging to infinity. Additionally, we impose restrictions upon
values of $t_j$ and functions $\varphi(z), \psi(z)$, such that we can scale the
contour $C$ by $q$ and $q^{-1}$ with respect to the point $z=0$
without crossing any poles. Under these conditions, the operator ${\cal D}$
formally becomes hermitian with respect to the taken inner product:
$ \langle \varphi,{\cal D}\psi \rangle=\langle {\cal D}\varphi,\psi \rangle. $
However, this property is not unique---the weight function
in the inner product can be multiplied by any elliptic function $\rho(z)$,
$\rho(qz)=\rho(z)$, with an accompanying change of the contour of integration.

In a trivial way,  $f(z)=1$ is an eigenfunction of  ${\cal D}$
with the eigenvalue $\lambda=0$ (actually, it solves simultaneously two
such equations---the second equation is obtained by permutation of $q$
and $p$). Evidently, the norm of this eigenfunction equals
to the elliptic hypergeometric function, $\|1\|^2=V(\underline{t})$.
This relation holds for $|p|,|q|<1$. If we change the integration variable
in the taken inner product $z=e^{2\pi i u/\omega_2}$, then, instead of $f(z)=1$,
we could have chosen as an $\lambda=0$ eigenfunction of $\mathcal{D}$
(where the operator $T$ is acting now as a shift, $Tv(u)=v(u+\omega_1)$)
any function $h(u)$ with the property $h(u+\omega_1)=h(u)$, but then the
normalization of this function would not be related to $V(\underline{t})$
in a simple way. For a special choice of this $h(u)$, we can obtain
$\|h\|^2=v^{mod}(\underline{g};{\bf\omega})$,
the modified elliptic hypergeometric function
for which we can take $|q|=1$. Equivalently, we could have changed the
inner product by replacing the standard elliptic gamma functions
by their modified version  and considering the
pair of equations $\mathcal{D}v(u)=0$ and its
$\omega_1\leftrightarrow \omega_2$ permuted partner.
Similar picture holds in the multivariable case considered in the end of this paper.

Because of these relations between $V(\underline{t})$ and the Calogero-Sutherland
type models, it
is natural to expect that elliptic hypergeometric functions will play a
major role in the solution of the standard eigenvalue problem for the
operator $\cal D$. In particular, we conjecture that the $E_7$ group of
symmetries of $V(\underline{t})$ can be lifted to $E_8$
at the level of unconstrained Hamiltonian ${\cal H}$ and that there is some
direct relation of this model with the elliptic Painlev\'e equation
(for this it would be desirable to understand an analogue of the
Painlev\'e-Calogero correspondence principle \cite{LO,man:painleve}
at the level of finite difference equations).

\section{Biorthogonal functions}

\subsection{Difference equation and three term recurrence relation}

For $n=0,1,\ldots,$ we define a sequence of functions \cite{spi:theta2}
$$
R_n(z;q,p)={}_{12}V_{11}\left(
\frac{t_3}{t_4};\frac{q}{t_0t_4},\frac{q}{t_1t_4},\frac{q}{t_2t_4}, t_3z,
\frac{t_3}{z},q^{-n}, \frac{Aq^{n-1}}{t_4};q,p\right),
$$
where $A=\prod_{m=0}^4t_m.$
They solve the elliptic hypergeometric equation rewritten in the form
$$
{\cal D}_\mu f(z)=0,\quad
{\cal D}_\mu=V_\mu(z)(T-1)+V_\mu(z^{-1})(T^{-1}-1)+\kappa_\mu,
$$
where $Tf(z)=f(qz)$ and
\begin{eqnarray*}
V_\mu(z)&=&
{ \theta\Bigl(\frac{pq\mu z}{t_4},\frac{pq^2z}{A\mu},\frac{t_4 z}{q};p\Bigr)}
\frac{\prod_{r=0}^4\theta(t_r z;p)}{\theta(z^2,qz^2;p)},  \\
\kappa_\mu&=&
{\theta\Bigl(\frac{A\mu}{qt_4},\mu^{-1};p\Bigr)}\prod_{r=0}^3
\theta\left(\frac{t_rt_4}{q};p\right),
\end{eqnarray*}
provided we quantize one of the parameters \fbox{$\mu=q^n$} (``the spectrum").
Equivalently, this equation can be rewritten as a generalized eigenvalue problem
$$
{\cal D}_\eta f(z)=\lambda {\cal D}_\xi f(z)
$$
with the spectral variable lying on the elliptic curve
$$
\lambda=\frac{\theta(\frac{\mu A\eta}{qt_4},\frac{\mu}{\eta};p)}
{\theta(\frac{\mu A\xi}{qt_4},\frac{\mu}{\xi};p)},
\quad
\xi,\eta\in {\mathbb C},\quad \xi\ne  \eta p^n,
\frac{qt_4}{A\eta}p^n,\; n\in\Z,
$$
where $\xi$ and $\eta$ are gauge parameters.
Out of this representation one obtains formal biorthogonality
$\langle T_n, R_m\rangle=0$ for $n\ne m$, where
$\langle \cdot, \cdot\rangle$ is some inner product
and $T_n(z;q,p)$ is a solution of a dual generalized eigenvalue problem.

{}From the elliptic hypergeometric equation one can derive also
the three-term recurrence relation
\begin{eqnarray*}\textstyle
&&(\gamma(z)-\alpha_{n+1})\rho(Aq^{n-1}/t_4)\Big(R_{n+1}(z;q,p)-R_{n}(z;q,p)\Big)\\
&& \makebox[2em]{}
+\ (\gamma(z)-\beta_{n-1})\rho(q^{-n})\Big(R_{n-1}(z;q,p)-R_{n}(z;q,p)\Big)\\
&&\makebox[4em]{}
+\ \delta\big(\gamma(z)-\gamma(t_3)\big)R_{n}(z;q,p)=0,
\end{eqnarray*}
with the initial conditions $R_{-1}=0,\; R_0=1$ and
\begin{eqnarray*}
&& \rho(x)=\frac{\theta(x,\frac{t_3}{t_4 x},\frac{qt_3}{t_4 x},
\frac{qx}{t_0t_1},\frac{qx}{t_0t_2},\frac{qx}{t_1t_2},
\frac{q^2\eta x}{A},\frac{q^2 x}{A\eta};p)}
{\theta(\frac{qt_4 x^2}{A},\frac{q^2 t_4x^2}{A};p)},\\
&& \makebox[1,4em]{}
\delta= \theta\left(\frac{q^2t_3}{A},\frac{q}{t_0t_4},\frac{q}{t_1t_4},
\frac{q}{t_2t_4},t_3\eta ,\frac{t_3}{\eta};p\right),
\\ &&
\gamma(z)=\frac{\theta(z\xi,z/\xi;p)}{\theta(z\eta,z/\eta;p)},\quad
\alpha_n=\gamma(q^n/t_4),\quad \beta_n=\gamma(q^{n-1}A).
\end{eqnarray*}
Since the whole $z$-dependence in this relation is concentrated in the
$\gamma(z)$ function, $R_n(z;q,p)$ are rational functions of $\gamma(z)$ with poles at
$\gamma(z)=\alpha_1,\ldots,\alpha_n$.

{}From the general theory of biorthogonal
rational functions \cite{zhe:bio} it follows that $R_n(z;q,p)$ can be orthogonal to
a rational function $T_n(z;q,p)$ with poles at $\gamma(z)=\beta_1,\ldots,\beta_n$.
The involution $t_4 \to pq/A$ permutes $\alpha_n$ and $\beta_n$, therefore the dual
functions are obtained after an application of this transformation to $R_n(z;q,p)$:
$$
T_n(z;q,p)=
{}_{12}V_{11}\left(\frac{At_3}{q};\frac{A}{t_0},\frac{A}{t_1},\frac{A}{t_2},
t_3z,\frac{t_3}{z},q^{-n},\frac{Aq^{n-1}}{t_4};q,p\right),
$$
where the $p$-dependence in parameters drops out due to the total ellipticity
property (in particular, we have $R_n(pz;q,p)=R_n(z;q,p)$).

\subsection{Two-index biorthogonality}

Let us denote the operator $\mathcal{D}_\mu$ introduced above as $\mathcal{D}_\mu^{q,p}$.
Then the product $R_{nm}(z)\equiv R_n(z;q,p)\cdot R_m(z;p,q)$
solves two generalized eigenvalue problems
$$
{\cal D}_\mu^{q,p}f(z)=0, \qquad {\cal D}_\mu^{p,q}f(z)=0
$$
with the spectrum \fbox{$\mu=q^np^m$.} Similar property is valid for the
dual product $T_{nm}(z)\equiv T_n(z;q,p)\cdot T_m(z;p,q)$ for a different choice
of parameters in $\mathcal{D}_\mu^{q,p}$.

\begin{thm}{\bf  (Two-index biorthogonality \cite{spi:theta2})} 
\par
If we denote
\begin{eqnarray*}
\Delta(z,\underbar{$t$})&=&\frac{(q;q)_\infty(p;p)_\infty}{4\pi i}\frac{\prod_{m=0}^4 {\mit\Gamma}(t_mz,t_mz^{-1})}
{{\mit\Gamma}(z^2,z^{-2}, Az, Az^{-1})},\\
{\cal N}(\underbar{$t$})&=&
\frac{\prod_{0\le m< k\le 4} {\mit\Gamma}(t_mt_k)}
   {\prod_{m=0}^4 {\mit\Gamma}(At_m^{-1})},
\end{eqnarray*}
where $|q|,|p|<1$, $|t_m|<1$, $|pq|<|A|$, then
\vskip 1mm \centerline{\fbox{$\displaystyle
\int_{C_{mn,kl}}T_{nl}(z)R_{mk}(z)\Delta(z,\underbar{$t$})\frac{dz}{z}
=h_{nl}\,{\cal N}(\underbar{$t$})\,\delta_{mn}\,\delta_{kl},
$}}\vskip 1mm\noindent
where $C_{mn,kl}$ is a contour separating points
$$t_jp^aq^b \,(j=0,1,2,3),\ t_4p^{a-k}q^{b-m}, p^{a+1-l}q^{b+1-n}/A,\quad
a,b\in {\mathbb N},
$$
from their $z\to z^{-1}$ reciprocals and normalization constants
\begin{eqnarray*}
h_{nl}&=&h_n(q,p)\cdot h_m(p,q),\\
h_n(q,p)&=&\frac{\theta(A/qt_4;p)
(q,qt_3/t_4,t_0t_1,t_0t_2,t_1t_2,At_3)_n\,q^{-n}}
{\theta(Aq^{2n}/t_4q;p) (1/t_3t_4,t_0t_3,t_1t_3,t_2t_3,A/qt_3,A/qt_4)_n}.
\end{eqnarray*}
\end{thm}
Only for $k=l=0$ there exists the $p\to 0$ limit and functions
$R_n(z;q,0)$ and $T_n(z;q,0)$ coincide with the Rahman's family of continuous
${}_{10}\varphi_9$ bi\-or\-tho\-go\-nal rational functions \cite{rah:integral}.
Note also that only for $k=l=0$ or $n=m=0$ we have rational functions of
some argument depending on $z$;
the general functions $R_{nm}(z)$ and $T_{nm}(z)$ should be
considered  as some meromorphic functions of $z$ with essential singularities
at $z=0$ and $z=\infty$.

For some quantized values of  $z$ and one of the parameters $t_j$
the functions $R_n(z;q,p)$ and $T_n(z;q,p)$ are reduced to the finite dimensional
set of biorthogonal rational functions constructed by Zhedanov and the author
in \cite{sz:cmp}. They generalize to the elliptic level Wilson's family of
discrete very well poised $_{10}\varphi_9$ biorthogonal functions \cite{wil:orthogonal}.
As described by Zhedanov at this workshop \cite{zhe:rims}, these functions have
found nice applications within the general Pad\'e interpolation scheme.

Functional solutions of the elliptic hypergeometric equation open the road to
construction of the associated biorthogonal functions following the
procedure described in \cite{ism-rah:associated} and this is one of the
interesting open problems for the future. A terminating continued fraction
generated by the three term recurrence relation described above has been calculated
in \cite{spi-zhe:theory}. It is expressed in terms of a terminating $_{12}V_{11}$
series and, again, the function $V(\underline{t})$ is expected to appear in the
description of non-terminating  convergent continued fractions generalizing
$q$-hypergeometric examples of \cite{gup-mas:contiguous}.

\subsection{The unit circle case}

In order to describe biorthogonal functions for which the measure
is defined by the modified elliptic beta integral,
we parametrize  $t_j= e^{2\pi i g_j/\omega_2}$ and introduce new notation
for the functions $R_n(z;q,p)$:
\begin{eqnarray*}
&& r_n(u;\omega_1,\omega_2,\omega_3)={}_{12}V_{11}
\Big(e^{2\pi i(g_3-g_4)/\omega_2};
e^{2\pi i(\omega_1-g_0-g_4)/\omega_2},
\\ && \makebox[2em]{} e^{2\pi i(\omega_1-g_1-g_4)/\omega_2},
e^{2\pi i(\omega_1-g_2-g_4)/\omega_2},
e^{2\pi i ({\cal A}+(n-1)\omega_1-g_4)/\omega_2},
\\ && \makebox[2em]{}
e^{-2\pi i n\omega_1/\omega_2},e^{2\pi i(g_3+u)/\omega_2},
e^{2\pi i(g_3-u)/\omega_2};e^{2\pi i \omega_1/\omega_2},
e^{2\pi i \omega_3/\omega_2}\Big),
\end{eqnarray*}
where ${\cal A}=\sum_{j=0}^4 g_j$. Similarly, we redenote the functions
$T_n(z;q,p)$ as $s_n(u; \omega_1,\omega_2,\omega_3)$.

The $q\leftrightarrow p$ symmetric situation (the standard set of biorthogonal
functions with $|p|, |q|<1$) is defined as the
$\omega_1\leftrightarrow\omega_3$ symmetric product of these functions:
$$
R_{nm}(e^{2\pi i u/\omega_2})=
r_n(u;\omega_1,\omega_2,\omega_3)\cdot r_m(u;\omega_3,\omega_2,\omega_1),
$$
with a similar relation for $T_{nm}(e^{2\pi i u/\omega_2})$.
As described above, we have the biorthogonality relations
$\langle T_{nl},R_{mk}\rangle= h_{nl}
\delta_{mn}\delta_{kl},$ where $\langle 1,1\rangle=1$ coincides with
the normalized standard elliptic beta integral with a special
contour of integration $C_{mn,kl}$ and
$$
h_{nl}=h_n(\omega_1,\omega_2,\omega_3)h_l(\omega_3,\omega_2,\omega_1)
$$
with $h_n(\omega_1,\omega_2,\omega_3)\equiv h_n(q,p)$. These functions are
modular invariant: $r_n(u;\omega_1,\omega_2,\omega_3)
=r_n(u;\omega_1,-\omega_3,\omega_2)$, $h_n(\omega_1,\omega_2,\omega_3)=
h_n(\omega_1,-\omega_3,\omega_2)$.

In the unit circle case we define functions
\begin{eqnarray*}
&& r_{nm}^{mod}(u)=r_n(u;\omega_1,\omega_2,\omega_3)\cdot
r_m(u;\omega_2,\omega_1,\omega_3),
\\ && s_{nm}^{mod}(u)=s_n(u;\omega_1,\omega_2,\omega_3)\cdot
s_m(u;\omega_2,\omega_1,\omega_3),
\end{eqnarray*}
which are now symmetric with respect to the permutations
 $\omega_2\leftrightarrow \omega_1$ and $n\leftrightarrow m$. These functions
satisfy the biorthogonality relations
$\langle s_{nl}^{mod},r_{mk}^{mod}\rangle =
h_{nl}^{mod}\delta_{mn}\delta_{kl}$,
where $\langle 1,1 \rangle=1$ coincides with the normalized modified elliptic
beta integral with the integration contour $\tilde{C}_{mn,kl}$ chosen
in an appropriate way and
$$
h_{nl}^{mod}=h_n(\omega_1,\omega_2,\omega_3)\cdot h_l(\omega_2,\omega_1,\omega_3).
$$

In sharp difference from the previous case, the limit $p\to 0$
(taken in such a way that simultaneously $r\to 0$, i.e.
Im$(\omega_3/\omega_1),$ Im$(\omega_3/\omega_2)\to +\infty$)
exists for all values of indices $n,l,k,m$ and we obtain:
\begin{eqnarray*}
&&r_{nm}(u;\omega_1,\omega_2)={}_{10}W_9\Big(e^{2\pi i(g_3-g_4)/\omega_2};\ldots,
e^{2\pi i(g_3-u)/\omega_2};q,q\Big)\\
&&\makebox[4em]{}\times
{}_{10}W_9\Big(e^{2\pi i(g_3-g_4)/\omega_1};\ldots,
e^{2\pi i(g_3-u)/\omega_1}; \tilde{q}^{-1}, \tilde{q}^{-1}\Big).
\end{eqnarray*}
Their partners from the dual space $s_{nm}(u;\omega_1,\omega_2)$ are defined
in a similar way. These functions $r_{nm}(u;\omega_1,\omega_2)$ and
$s_{nm}(u;\omega_1,\omega_2)$ are not rational functions of some
particular combination of the variable $u$ for $n, m\neq 0$.
They satisfy the two-index biorthogonality relations
$$
\langle r_{nl},s_{mk}\rangle =\nu_{nl} \delta_{mn}\delta_{kl},
$$
where $\nu_{nl}$ are obtained from $h_{nl}^{mod}$ after setting
Im$(\omega_3/\omega_1),$ Im$(\omega_3/\omega_2)\to +\infty$
and $\langle 1,1\rangle=1$ coincides with the normalized ``unit circle" partner of the
Rahman's integral \cite{sto:hyperbolic} with a special contour of integration.
Further simplification of these relations to the Askey-Wilson polynomials
level is highly non-trivial due to some problems with the
convergence of the integral and requires a thorough investigation.
In a similar way it is possible to define unit circle partners of the
Rains' multivariable generalization of the author's univariate biorthogonal functions
\cite{rai:abelian} as well as their limiting
two-index $q$-biorthogonal functions.

\section{Multiple elliptic beta integrals}

\subsection{General definition}

Multiple integrals
$$
\int_{D}\Delta(u_1,\ldots, u_n)\, du_1\cdots du_n,
$$
where $D\subset {\mathbb C}^n$ are some $n$-dimensional cycles,
are called elliptic hypergeometric integrals if $\Delta(u_1,\ldots,u_n)$
are meromorphic functions of $u_1,\ldots, u_n$ satisfying the
following system of equations
$$
\Delta(u_1,\ldots, u_k+\omega_1,\ldots  ,u_n)=h^{(k)}(u_1,\ldots, u_n)\,
\Delta(u_1,\ldots, u_n),
$$
where $h^{(k)}(\underline{u})$, $k=1,\ldots  ,n,$ are elliptic functions of all
$u_j$, i.e.,
$$
h^{(k)}(u_j+\omega_2)=h^{(k)}(u_j+\omega_3)=h^{(k)}(u),\quad
{\rm Im}(\omega_2/\omega_3)\ne 0.
$$

This is a ``broad'' definition of the integrals introduced in \cite{spi:theta2};
one can make it ``narrow'' by tripling the number of equations for
$\Delta(\underline{u})$ using the shifts by all quasiperiods $\omega_i$.

In order to describe general possible forms of the integrand,
we need an elliptic extension of the Ore-Sato theorem on the
general form of terms in plain hypergeometric series (see, e.g.,
\cite{ggr:general}).
For all ``good" known elliptic hypergeometric integrals, the kernels
$\Delta(\underline{u})$ are equal to ratios of elliptic gamma functions $\Gamma(z;q,p)$
with an integer power dependence on the variables $z_j=e^{2\pi iu_j/\omega_2}$.
However, in general case we can multiply the integrands
by elliptic functions of all $u_j$'s with the periods
$\omega_2, \omega_3$ which do not have such a representation.

Multiple elliptic hypergeometric series are defined in a similar way.
It is simply necessary to replace integrals by discrete sums over
some sublattices of $u_1,\ldots,u_n\in\Z^n$ keeping other properties
of $\Delta(\underline{u})$. We shall not consider them in the present review.

The most interesting elliptic hypergeometric integrals are related to
multiple generalizations of the elliptic beta integral, which are split
formally into three different groups. Type I integrals contain $2n+3$ free parameters
and bases $p$ and $q$ and their proofs use in one or another way
analytical continuation procedure over discrete values of parameters.
Type II integrals contain less than $2n+3$ free parameters and
they can be proved by purely algebraic means on the basis of type I integrals.
Finally, type III elliptic beta integrals arise through
computations of $n$-dimensional determinants with entries composed of
one-dimensional integrals. It goes without saying that all these
integrals have their partners expressed in terms of the modified
elliptic gamma function.

\subsection{Integrals for the root system $C_n$}

In order to define $n$-di\-men\-sio\-nal type I elliptic beta integral
for the root system $C_n$ (abbreviated as the $C_I$ integral),
we take bases $|p|,|q|<1$ and parameters
$t_1,\ldots ,$ $t_{2n+4}\in\C$ such that
$\prod_{j=1}^{2n+4}t_j=pq$ and $|t_1|,\ldots,|t_{2n+4}|<1$.
\begin{thm}{\bf (Type I $C_n$ elliptic beta integral \cite{die-spi:selberg})}
\vskip 2mm \centerline{\fbox{$\displaystyle
\kappa_n^C\int_{\T^n}
\prod_{j=1}^n\frac{\prod_{i=1}^{2n+4}\Gamma(t_iz_j^{\pm})}
{\Gamma(z_j^{\pm 2})}
\prod_{1\leq i<j\leq n} \frac{1}{\Gamma(z_i^{\pm}z_j^{\pm})}
\frac{dz}{z}
=\prod_{1\leq i<j\leq 2n+4}\Gamma(t_it_j),
$}}\vskip 1mm \noindent
where $\Gamma(z)\equiv \Gamma(z;q,p)$ and
\vskip 1mm \centerline{\fbox{$\displaystyle
\kappa_n^{C}=\frac{(p;p)_\infty^n(q;q)_\infty^n}{(2\pi i)^n2^n n!}.
$}}
\end{thm}
Different complete proofs of this formula were given by Rains \cite{rai:trans}
and the author \cite{spi:short}.
In the limit $p\to 0$ it is reduced to one of the Gustafson results
\cite{gus:some1}. Its modified elliptic gamma function partner has
been established by the author \cite{spi:short} together with
its $q$-degeneration valid for $|q|\leq 1$ (which we skip for
brevity).

Type II integral for this root system (abbreviated as the $C_{II}$ integral)
depends on seven parameters $t$ and
$t_m$, $m=1,\ldots ,6,$  and bases $q,p$ constrained by one relation.
It can be derived as a consequence of the $C_I$ integral.

\begin{thm}{\bf (Type II $C_n$ elliptic beta integral \cite{die-spi:elliptic})}
\par
Let nine complex parameters $t, t_m (m=1,\ldots , 6), p$ and $q$ be constrained
by the conditions $|p|, |q|,$ $|t|,$ $|t_m| <1,$ and
$t^{2n-2}\prod_{m=1}^6t_m=pq$. Then,
\vskip 1mm \centerline{\fbox{$
\begin{array}{c}
\displaystyle
\kappa_n^C\int_{\T^n} \prod_{1\leq j<k\leq n}
\frac{\Gamma(tz_j^\pm z_k^\pm)}{\Gamma(z_j^\pm z_k^\pm)}
\prod_{j=1}^n\frac{\prod_{m=1}^6\Gamma(t_mz_j^\pm)}{\Gamma(z_j^{\pm2})}
\frac{dz}{z}
\\ \displaystyle
= \prod_{j=1}^n\left(\frac{\Gamma(t^j)}{\Gamma(t)}
\prod_{1\leq m<s\leq 6}\Gamma(t^{j-1}t_mt_s)\right).
\end{array}
$}} \vskip 1mm
\end{thm}

This is an elliptic analogue of the Selberg integral which appears after
a number of reductions, the first step being the $p\to0$ limit
leading to one of the Gustafson's integrals \cite{gus:some2}.
In order to take this limit it is necessary to express $t_6$ in terms of
other parameters and remove the multipliers $pq$ by the
inversion formula for $\Gamma(z;q,p)$ (see \cite{die-spi:selberg}). During this
procedure the integral takes a less symmetric form---in the given form
it has the explicit $S_6$ symmetry in parameters (see \cite{rai:trans}).
For the modified version of this integration formula valid for $|q|\leq 1$
and its $q$-degeneration, see \cite{die-spi:unit}.

Presently the author knows only one type III elliptic beta integral \cite{spi:theta2}.
It is ascribed to the $C_n$  root system (we abbreviate it as the $C_{III}$
integral) and it is computed by evaluation
of a determinant of the univariate elliptic beta integrals which is reduced
to the computation of the Warnaar's determinant \cite{war:summation}.
We skip it for brevity, but it is expected that there are much more such integrals
due to the universality of the method used for their derivation (see, e.g., \cite{tv})
and existence of several nice exact determinant evaluations for elliptic theta
functions.

\subsection{Integrals for the root system $A_n$}

Classification of the $A_n$ elliptic beta integrals follows the
same line as in the $C_n$ case. We start from the description
of the simplest type I integral introduced by the author in
\cite{spi:theta2}, which we symbolize as $A_I^{(1)}$.

\begin{thm}{\bf (The $A_I^{(1)}$ integral \cite{spi:theta2})}
\vskip 1mm \centerline{\fbox{$
\begin{array}{c}
\displaystyle
\kappa_n^A\int_{{\mathbb T}^n}
\prod_{1\le j<k\le n+1}\frac{1}{{\mit\Gamma}(z_iz_j^{-1},z_i^{-1}z_j)}
\,\prod_{j=1}^{n+1}\prod_{m=1}^{n+2}{\mit\Gamma}(s_mz_j,t_mz_j^{-1})
\,\frac{dz}{z}
\\ \displaystyle
=\prod_{m=1}^{n+2} {\mit\Gamma}(Ss_m^{-1},Tt_m^{-1})
\prod_{k,m=1}^{n+2} {\mit\Gamma}(s_kt_m),
\end{array}
$}} \vskip 1mm
\noindent
where \fbox{$z_1z_2\cdots z_{n+1}=1$} and
\vskip 1mm \centerline{\fbox{$\displaystyle
\kappa_n^A=\frac{(p;p)_\infty^n(q;q)_\infty^n}{(2\pi i)^n(n+1)!}
$}} \vskip 1mm\noindent
with the parameters satisfying the constraints $|t_m|, |s_m|<1$, $m=1,\ldots, n+2,$
and \fbox{$ST=pq$,} $S=\prod_{m=1}^{n+2}s_m$, $T=\prod_{m=1}^{n+2}t_m$.
\end{thm}

For complete proofs of this formula, see \cite{rai:trans,spi:short}.
Here we have a split of $2n+4$ parameters (homogeneous in the $C_n$ case)
with one constraint into two homogeneous groups with $n+2$ entries in each group.
The $p\to 0$ limiting value of this integral was derived by Gustafson \cite{gus:some1}.
The unit circle analogue together with the appropriate $q$-degeneration valid for
$|q|\leq 1$ were derived in \cite{spi:short}. Another type I $A_n$ integral
is described below.

There are several type II integrals on the $A_n$ root system,
the first of which we abbreviate as $A_{II}^{(1)}$.
For its description we define the kernel
\vskip 1mm \centerline{\fbox{$\displaystyle
\Delta^{(1)}_{II}(\underline{z})=\prod_{1\leq i<j\leq n+1}
\frac{\Gamma(tz_iz_j)}{\Gamma(z_iz_j^{-1},z_i^{-1}z_j) }
 \prod_{j=1}^{n+1} \prod_{k=1}^{n+1}\Gamma(t_kz_j^{-1})
\prod_{i=1}^4\Gamma(s_iz_j),
$}} \vskip 1mm\noindent
where $t^{n-1}\prod_{k=1}^{n+1}t_k\prod_{i=1}^4s_i=pq$ and $\prod_{j=1}^{n+1}z_j=1$.

\begin{thm} {\bf (The $A_{II}^{(1)}$ integral \cite{spi:theta2})}
\par
As a consequence of the $C_I$ and $A_I^{(1)}$ integration formulas, we have
for odd $n$
\vskip 1mm \centerline{\fbox{$
\begin{array}{c}
\displaystyle
 \kappa_n^A\int_{\mathbb{T}^n}\Delta^{(1)}_{II}(\underline{z})\frac{dz}{z}
= \frac{\Gamma(t^{\frac{n+1}{2}},A)}{\Gamma(t^{\frac{n+1}{2}}A)}
\prod_{k=1}^{n+1}\prod_{i=1}^4\Gamma(t_ks_i)
 \\  \displaystyle  \makebox[4em]{} \times
\prod_{1\leq j<k\leq n+1} \Gamma(tt_jt_k)
\prod_{1\leq i<m\leq 4}\Gamma(t^{\frac{n-1}{2}}s_is_m).
\end{array}
$}} \vskip 1mm
\noindent
where $A=\prod_{k=1}^{n+1}t_k$.

 For even $n$, we have
\vskip 1mm \centerline{\fbox{$
\begin{array}{c}
\displaystyle
\kappa_n^A\int_{\mathbb{T}^n}\Delta^{(1)}_{II}(\underline{z})\frac{dz}{z}
= \Gamma(A )\prod_{k=1}^{n+1}\prod_{i=1}^4\Gamma(t_ks_i)
\\ \displaystyle \makebox[1em]{} \times
\prod_{1\leq j<k\leq n+1} \Gamma(tt_jt_k)
\prod_{i=1}^4\frac{\Gamma(t^{\frac{n}{2}}s_i)}
{\Gamma(t^{\frac{n}{2}}As_i) }.
\end{array}
$}} \vskip 1mm
\noindent
\end{thm}

These formulas contain only $n+5$ free parameters.
In the  $p\to 0$ limit they are reduced to the main result of
\cite{gus-rak:beta}.

We abbreviate the second type II $A_n$ integral  as  $A_{II}^{(2)}$. For
its description we need the kernel
\vskip 1mm \centerline{\fbox{$\displaystyle
\Delta_{II}^{(2)}(\underline{z})= \prod_{1\leq i<j\leq n+1}
\frac{\Gamma(tz_iz_j,sz_i^{-1}z_j^{-1})}
{\Gamma(z_iz_j^{-1},z_i^{-1}z_j)}
\prod_{j=1}^{n+1}\prod_{k=1}^3\Gamma(t_kz_j,s_kz_j^{-1}),
$}} \vskip 1mm\noindent
where ten variables $p,q, t,s, t_1,t_2,t_3,s_1,s_2, s_3\in\C$ satisfy one
constraint $(ts)^{n-1}\prod_{k=1}^3t_ks_k=pq$.

\begin{thm} {\bf (The $A_{II}^{(2)}$ integral \cite{spi:theta2})}
\par
As a consequence of the $A_I^{(1)}$, $C_I,$ and $C_{II}$ integration formulas,
we have an additional type II elliptic beta integral for the $A_n$ root system.
For odd $n$, we have
\vskip 1mm \centerline{\fbox{$
\begin{array}{c}
\displaystyle
\kappa_n^A\int_{\mathbb{T}^n}\Delta_{II}^{(2)}(\underline{z}) \frac{dz}{z}
=\Gamma(t^{\frac{n+1}{2}},s^{\frac{n+1}{2}})
\prod_{1\leq i<k\leq 3} \Gamma(t^{\frac{n-1}{2}}t_it_k,s^{\frac{n-1}{2}}s_is_k)
\\ \displaystyle \makebox[4em]{} \times
\prod_{j=1}^{(n+1)/2}\prod_{i,k=1}^3\Gamma((ts)^{j-1}t_is_k)
\\  \displaystyle \makebox[4em]{} \times
\prod_{j=1}^{(n-1)/2}\left(\Gamma((ts)^j)\prod_{1\leq i<k\leq 3}
\Gamma(t^{j-1}s^jt_it_k,t^js^{j-1}s_is_k)\right).
\end{array}
$}} \vskip 1mm
\noindent
For even $n$, we have
\vskip 1mm \centerline{\fbox{$
\begin{array}{c}
\displaystyle
 \kappa_n^A\int_{\mathbb{T}^n}\Delta_{II}^{(2)}(\underline{z}) \frac{dz}{z}
= \prod_{i=1}^3\Gamma(t^{\frac{n}{2}}t_i,s^{\frac{n}{2}}s_i) \qquad\qquad\qquad
\\ \displaystyle \makebox[4em]{} \times
\Gamma(t^{\frac{n}{2}-1}t_1t_2t_3,s^{\frac{n}{2}-1}s_1s_2s_3)
\prod_{j=1}^{n/2}\Biggl(\Gamma((ts)^j)
\\ \displaystyle 
\times \prod_{i,k=1}^3\Gamma((ts)^{j-1}t_is_k)
\prod_{1\leq i<k\leq 3}\Gamma(t^{j-1}s^jt_it_k,t^js^{j-1}s_is_k)\Biggr).
\end{array}
$}} \vskip 1mm
\noindent
In this and previous theorems we assume constraints on the parameters
guaranteeing that all sequences of integrands' poles converging
to zero (or their reciprocals) lie within (or outside) of $\T$.
\end{thm}

This theorem formulas contain only seven free parameters.
In the $p\to 0$ limit we obtain one of the integrals in \cite{gus:some2}.

Recently, Warnaar and the author have found a complementary type I elliptic
beta integral for the $A_n$ root system.
\begin{thm} {\bf (The $A_{I}^{(2)}$ integral \cite{spi-war:inversions})}
\vskip 1mm \centerline{\fbox{$
\begin{array}{c}
\displaystyle
\kappa_n^A\int_{{\mathbb T}^n}
\prod_{1\le i<j\le n+1}\frac{{\mit\Gamma}(Sz_i^{-1}z_j^{-1})}
{{\mit\Gamma}(z_iz_j^{-1},z_i^{-1}z_j)}
\prod_{j=1}^{n+1}
\frac{\prod_{k=1}^{n}{\mit\Gamma}(t_kz_j)\prod_{m=1}^{n+3}{\mit\Gamma}(s_mz_j^{-1})}
{\prod_{k=1}^{n}{\mit\Gamma}(St_kz_j^{-1})} \frac{dz}{z}
\\ \displaystyle
= \prod_{k=1}^{n} \prod_{m=1}^{n+3}
\frac{{\mit\Gamma}(t_ks_m)}{{\mit\Gamma}(St_ks_m^{-1})}
\prod_{1\le l<m\le n+3} {\mit\Gamma}(Ss_l^{-1}s_m^{-1}),  \makebox[2em]{}
\end{array}
$}} \vskip 1mm
\noindent
where $|t_k|<1\ (k=1,2,\ldots ,n)$, $|s_m|<1\ (m=1,2,\ldots,n+3)$,
$|pq|<|t_j S|$, $S=\prod_{m=1}^{n+3}s_m$,
and $z_1\cdots z_{n+1}=1$.
\end{thm}
Here we have a split of $2n+3$ independent parameters into two groups
with $n$ and $n+3$ homogeneous entries. 
This integration formula appeared to be new even in the $p\to0$ limit
as well as in its further degeneration to the plain hypergeometric level $q\to 1$.
Its unit circle analogue valid for $|q|\leq 1$ is constructed in \cite{spi:short}.

For each of the described integrals we can apply the residue calculus similar
to the one described above in the univariate case and derive
summation formulas for particular multiple elliptic
hypergeometric series on root systems generalizing the Frenkel-Turaev sum.
For the $C_I$ integral, the corresponding
formula was derived by van Diejen and the author \cite{die-spi:modular}
and its recursive proof was given by Rosengren \cite{ros:elliptic}.
For the $C_{II}$ integral, the corresponding sum was conjectured first by Warnaar
\cite{war:summation}, it was deduced from the residue calculus by van Diejen and
the author \cite{die-spi:elliptic} and proven recursively by Rosengren
\cite{ros:proof}. The $A_I^{(1)}$ resides sum was deduced by the author \cite{spi:theta2},
leading to an elliptic generalization of the Milne's sum \cite{mil:multiple}.
Residue calculus for the $A_I^{(2)}$ integral performed by Warnaar and the
author \cite{spi-war:inversions} leads to an elliptic generalization of
the Bhatnagar-Schlosser ``$D_n$" summation formula \cite{BS}. These elliptic
$A_I^{(1)}$ and $A_I^{(2)}$ summation formulas were proven first
inductively by Rosengren \cite{ros:elliptic}. A summation formula following
from the $A_{II}^{(1)}$ integral was conjectured by the author \cite{spi:theta2},
but it still remains unproven. Residue calculus for the $C_{III}$
integral is expected to lead to a Warnaar's sum \cite{war:summation}, but
this question was not investigated either.

All the described integrals are expected to serve as measures in the orthogonality
relations for some biorthogonal functions. A program of searching
multivariable analogues of the $_{12}V_{11}$ biorthogonal functions was put forward in
\cite{die-spi:elliptic,spi:theta1}. The first example of a multivariable
extension of the author's two-index continuous biorthogonal
functions was found by Rains \cite{rai:abelian} on the basis of the
$C_{II}$ elliptic beta integral (these functions generalize also the
Okounkov's interpolating polynomials \cite{oku}).

The notion of root systems provides the main guiding principle in the
construction of multiple elliptic beta integrals. Although this connection is not
straightforward, it is natural to expect that there exist other such
integrals attached, in particular, to the exceptional Lie algebras.
In this respect it is worth analyzing whether all multiple Askey-Wilson
type integrals classified by Ito \cite{ito} admit a further lift up to the
levels of Rahman's $q$-beta integral and the author's elliptic beta integral.

\section{Univariate integral  Bailey chains}

The Bailey chains techniques is well known as a powerful tool
for derivation of infinite sequences of identities for series
of hypergeometric type \cite{aar:special}. The most general known
$q$-hypergeometric Bailey chain was proposed by Andrews \cite{and:bailey}.
It is related to the Bressoud's matrix inverse \cite{bre} and has at the
bottom the original constructions by Rogers and Bailey used for proving
the Rogers-Ramanujan identities \cite{bm}. It was generalized
to the elliptic hypergeometric series by the author
\cite{spi:bailey1} (for some further developments in this direction,
see \cite{war:extensions}).
We shall not describe these chains here, although they have quite interesting
consequences. Instead, we present Bailey chains for integrals
discovered in \cite{spi:bailey2}.

\underline{D{\scriptsize EFINITION}.} Two functions $\alpha(z,t)$ and
 $\beta(z,t)$ form an elliptic integral Bailey pair with respect to
the parameter $t$, if
\vskip 1mm \centerline{\fbox{$\displaystyle
\beta(w,t)=\kappa\int_{{\mathbb T}}{\mit\Gamma}(tw^{\pm}z^{\pm})\alpha(z,t)\frac{dz}{z}.
$}}\vskip 1mm

\begin{thm} {\bf (First integral Bailey lemma \cite{spi:bailey2})}
\par
For a given integral Bailey pair $\alpha(z,t)$, $\beta(z,t)$ with respect to $t$,
the functions
\begin{eqnarray*}
&& \alpha'(w,st)=\frac{{\mit\Gamma}(tuw^{\pm})}
{{\mit\Gamma}(ts^2uw^{\pm})}\alpha(w,t),\\
&& \beta'(w,st)=\kappa \frac{{\mit\Gamma}(t^2s^2,t^2suw^{\pm})}
{{\mit\Gamma}(s^2,t^2,suw^{\pm})}
\int_{{\mathbb T}} \frac{{\mit\Gamma}(sw^{\pm}x^{\pm},ux^{\pm})}
{{\mit\Gamma}(x^{\pm 2},t^2s^2ux^{\pm})}\beta(x,t)\frac{dx}{x},
\end{eqnarray*}
where $w\in {\mathbb T}$, form a new Bailey pair with respect to
the parameter $st$.
\end{thm}
The proof is quite simple, it is necessary to substitute the key relation
for $\beta(x,t)$ into the definition of $\beta'(w,st)$, to change the
order of integrations, and to apply the
elliptic beta integral (under some mild restrictions upon parameters).
Note that these substitution rules introduce two new parameters $u$ and $s$
into the Bailey pairs at each step of their iterative application.

\begin{thm}{\bf (Second integral Bailey lemma \cite{spi:bailey2})}

For a given integral Bailey pair $\alpha(z,t)$, $\beta(z,t)$ with respect
to the parameter $t$, the functions
\begin{eqnarray*}
&& \alpha'(w,t)=\kappa \frac{{\mit\Gamma}(s^2t^2,uw^{\pm})}
{{\mit\Gamma}(s^2,t^2,w^{\pm 2},t^2s^2uw^{\pm})}
\int_{{\mathbb T}}
\frac{{\mit\Gamma}(t^2sux^{\pm},sw^{\pm}x^{\pm})}
{{\mit\Gamma}(sux^{\pm})}\alpha(x,st)\frac{dx}{x},\\
&& \beta'(w,t)=\frac{{\mit\Gamma}(tuw^{\pm})}{{\mit\Gamma}(ts^2uw^{\pm})}\beta(w,st)
\end{eqnarray*}
form a new Bailey pair with respect to $t$.
\end{thm}

It appears that these two lemmas are related to each other by inversion
of the integral operator figuring in the definition of integral
Bailey pairs \cite{spi-war:inversions}.
Application of these lemmas is algorithmic: one should take the initial
$\alpha(z,t)$ and $\beta(z,t)$ defined by the elliptic beta integral and
apply to them described transformations in all possible ways, which yields
a binary tree of identities for multiple elliptic hypergeometric integrals
of different dimensions. In particular, the very first step yields the
key transformation (i) for the elliptic hypergeometric function $V(\underline{t})$.
The residue calculus is supposed to recover elliptic Bailey chains for the
$_{r+1}V_r$ series \cite{spi:bailey1}.
We can take the limit $p\to 0$ and reduce all elliptic results 
to the level of standard $q$-hypergeometric integrals which admit further
simplification down to identities generated by the plain hypergeometric
beta integrals.

As to the unit circle case, we can start from the relation
$$
\tilde\beta(v,g)=\tilde\kappa\int_{-\omega_3/2}^{\omega_3/2}
G(g\pm v \pm u;{\bf \omega})\, \tilde \alpha(u,g)\frac{du}{\omega_2}
$$
and apply the modified elliptic beta integral for building needed analogues
of the Bailey lemmas. In this case, the $p,r\to 0$ limit brings in
identities for $q$-hypergeometric integrals defined over the non-compact
contour $\L$ with the kernels well defined for $|q|=1$.

\section{Elliptic Fourier-Bailey type integral transformations on root systems}

Similar to the situation with elliptic beta integrals,
the univariate integral transformation of the previous section has been
generalized by Warnaar and the author to root systems \cite{spi-war:inversions}.
It appears that in the multivariable setting the original space of functions and
its image can belong to different root systems.

For the $(A,A)$ pair of root systems, we take the space of meromorphic functions
$f_A(\underline{z};t)$ with $A_n$ symmetry in its variables $z_1,\ldots, z_{n+1}$,
$\prod_{j=1}^{n+1}z_j=1$, and define its image space by setting
\vskip 1mm \centerline{\fbox{$\displaystyle
\widehat{f}_A(\underline{w};t)=\kappa_n^A\int_D\rho(\underline{z},
\underline{w};t^{-1}) f_A(\underline{z};t)\frac{dz}{z},
$}}\vskip 1mm\noindent
where the kernel has the form
\vskip 1mm
\centerline{\fbox{$\displaystyle
\rho(\underline{z},\underline{w};t)=\frac{\prod_{i,j=1}^{n+1}{\mit\Gamma}(tw_i^{-1}z_j^{-1})}
{{\mit\Gamma}(t^{n+1})\prod_{1\le i<j\le n+1}{\mit\Gamma}(z_iz_j^{-1},z_i^{-1}z_j)}.
$}}\vskip 1mm\noindent
In a relatively general situation this map can be inverted explicitly.
\begin{thm} {\bf (The $(A,A)$ transform inversion \cite{spi-war:inversions})}
\par
For a suitable $n$-dimensional cycle $D$, the inverse of the $(A,A)$ transform
is given by the map
\vskip 1mm \centerline{\fbox{$\displaystyle
f_A(\underline{x};t)=\kappa_n^A \int_{{\mathbb T}^n}
\rho(\underline{w}^{-1},\underline{x}^{-1};t)\widehat{f}_A(\underline{w};t)\frac{dw}{w},
$}}\vskip 1mm\noindent
where it is assumed that functions $f_A(\underline{x};t)$ are analytical in a
sufficiently wide annulus encircling $\T$.
\end{thm}
The proof consists in a quite tedious  residue calculus
with the use of the $A_I^{(1)}$ integration formula.

In the $(A,C)$-case, we map functions $f_A(\underline{z};t)$
to its image space belonging to the $C_n$ root system:
\vskip 1mm \centerline{\fbox{$\displaystyle
\widehat{f}_C(\underline{w};t)=\kappa_n^A\int_D
\delta_A(\underline{z},\underline{w};t^{-1})f_A(\underline{z};t)\frac{dz}{z},
$}}\vskip 1mm\noindent
where the kernel has the form
\vskip 1mm \centerline{\fbox{$\displaystyle
\delta_A(\underline{z},\underline{w};t)=\frac{\prod_{i=1}^{n}
\prod_{j=1}^{n+1}{\mit\Gamma}(tw_i^{\pm}z_j^{\pm})}{\prod_{1\le i<j\le n+1}
{\mit\Gamma}(z_iz_j^{-1},z_i^{-1}z_j,t^{-2}z_iz_j,t^2z_i^{-1}z_j^{-1})}.
$}}\vskip 1mm\noindent
\begin{thm}{\bf (The $(A,C)$ transform inversion \cite{spi-war:inversions})}
\par
For a suitable $n$-dimensional cycle $D$, the inverse of the $(A,C)$
integral transform looks as follows
\vskip 1mm \centerline{\fbox{$\displaystyle
f_A(\underline{x};t)=\kappa_n^C \int_{{\mathbb T}^n}\delta_C(\underline{w},\underline{x};t)
\widehat{f}_C(\underline{w};t)\frac{dw}{w},
$}}\vskip 1mm\noindent
with the kernel
\vskip 1mm \centerline{\fbox{$\displaystyle
\delta_C(\underline{w},\underline{x};t)=
\frac{\prod_{i=1}^{n}\prod_{j=1}^{n+1}{\mit\Gamma}(tw_i^{\pm}x_j)}
{\prod_{i=1}^{n}{\mit\Gamma}(w_i^{\pm 2})
\prod_{1\le i<j\le n}{\mit\Gamma}(w_i^{\pm}w_j^{\pm})},
$}}\vskip 1mm\noindent
where it is assumed that functions $f_A(\underline{x};t)$ are analytical in a
sufficiently wide annulus containing $\T$.
\end{thm}
\begin{cor}
If we choose $\widehat{f}_C(\underline{w};t)$ such that the product
$\delta_C(\underline{w},\underline{x};t)\cdot\widehat{f}_C(\underline{w};t)$
is equal to the $C_I$ elliptic beta integral kernel, then the original relation
$\widehat{f}_C\sim \int_{D} \delta_A\cdot f_A\, dz/z$ defines the
$A_I^{(2)}$ integration formula.
\end{cor}

There are more such Fourier-Bailey type integral transforms with
explicit inversions some of which still are in the conjectural form.
All of them can be put into the integral Bailey chains setting
yielding many infinite sequences of transformations for the
elliptic hypergeometric integrals on root systems.

\section{Applications to the Calogero-Sutherland type models}

After discussing multiple elliptic beta integrals, we would like to return
to applications of elliptic hypergeometric functions to the
Ca\-lo\-ge\-ro-Su\-ther\-land type models \cite{spi:thesis}.

First, we define the inner product
$$
\langle \varphi, \psi \rangle^{I,II}=
\kappa_n^C
\int_{{\mathbb T}^n}\Delta^{I,II}(\underline{z},\underline{t})
\varphi(\underline{z})\psi(\underline{z})\,\frac{dz}{z}.
$$
Let us take the Hamiltonian of the van Diejen model \cite{die:integrability}
with the restriction $t^{2n-2}\prod_{m=1}^8t_m=p^2q^2$
\begin{eqnarray*}
&& {\cal D}_{II}=\sum_{j=1}^n\Big(A_j(\underline{z})(T_j-1)
+A_j(\underline{z}^{-1})(T_j^{-1}-1)\Big),
\\&&
 A_j(\underline{z})=\frac{\prod_{m=1}^8\theta(t_mz_j;p)}{\theta(z_j^2,qz_j^2;p)}
\prod_{k=1\atop \,\ \ne j}^n
\frac{\theta(tz_jz_k,tz_jz_k^{-1};p)}{\theta(z_jz_k,z_jz_k^{-1};p)}.
\end{eqnarray*}
Under some relatively mild restrictions upon parameters, this operator
is formally hermitian with respect to the above inner product,
$\langle \varphi, {\cal D}_{II}\psi \rangle^{II}
=\langle {\cal D}_{II}\varphi, \psi \rangle^{II},$
for the weight function
$$
\Delta^{II}(\underline{z},\underline{t})=\prod_{1\le j<k\le n}\!
\frac{{\mit\Gamma}(tz_j^{\pm}z_k^{\pm})}
{{\mit\Gamma}(z_j^{\pm}z_k^{\pm})}
\ \prod_{j=1}^n \frac{\prod_{k=1}^8{\mit\Gamma}(t_kz_j^{\pm})}
{{\mit\Gamma}(z_j^{\pm 2})}.
$$
Evidently, $f(\underline{z})=1$ is a $\lambda=0$ solution
of the standard eigenvalue problem
${\cal D}_{II}f(\underline{z})=\lambda f(\underline{z}).$
The norm of this eigenfunction
\vskip 1mm \centerline{\fbox{$\displaystyle
\|1\|^2=V(\underline{t};C_{II})=\kappa_n^C\int_{{\mathbb T}^n}
\Delta^{II}(\underline{z},\underline{t})\,\frac{dz}{z}
$}}\vskip 1mm\noindent
is a multivariable analogue of the elliptic hypergeometric
function $V(\underline{t})$ for the type II $C_n$ elliptic beta integral.

We conjecture that with all multiple elliptic beta integrals
one can associate Calogero-Sutherland type models in the described fashion.
Let us take the  weight function
$$
\Delta^{I}(\underline{z},\underline{t})=
\frac{1}{\prod_{1\le i<j\le n}{\mit\Gamma}(z_i^{\pm}z_j^{\pm})}
\prod_{j=1}^{n}\frac{\prod_{k=1}^{2n+6}{\mit\Gamma}(t_kz_j^{\pm})}
{{\mit\Gamma}(z_j^{\pm 2})}.
$$
We associate with it the Hamiltonian
\begin{eqnarray*}
&&\makebox[2em]{}
{\cal D}_I=\sum_{j=1}^n\Big(A_j(z)(T_j-1)+A_j(z^{-1})(T_j^{-1}-1)\Big),
\\&&
A_j(z)=\frac{\prod_{k=1}^{2n+6}\theta(t_kz_j;p)}{\theta(z_j^2,qz_j^2;p)}
\prod_{k=1\atop \,\ \ne j}^n
\frac{1}{\theta(z_jz_k^{\pm};p)},
\quad
\prod_{k=1}^{2n+6}t_k=p^2q^2,
\end{eqnarray*}
which is formally hermitian  with respect to the taken inner product,
$\langle \varphi, {\cal D}_I\psi \rangle^I=\langle {\cal D}_I\varphi, \psi \rangle^I,$
for some mild restrictions upon the parameters.
Again, $f(\underline{z})=1$ is a $\lambda=0$ eigenfunction of the operator
 ${\cal D}_I$ and its normalization
\vskip 1mm \centerline{\fbox{$\displaystyle
\|1\|^2=V(\underline{t};C_I)=\kappa_n^C\int_{{\mathbb T}^n}
\Delta^I(\underline{z},\underline{t})\,\frac{dz}{z}
$}}\vskip 1mm\noindent
defines type I generalization of the elliptic hypergeometric function
for the root system $C_n$. The functions $V(\underline{t};C_{I,II})$
were considered first by Rains \cite{rai:trans} in the context of symmetry
transformations for multiple elliptic hypergeometric integrals.
It is not difficult to define
their unit circle analogues which also play similar role in the context of
Calogero-Sutherland type models.

One can construct analogues of the $V(\underline{t})$ function for
multiple elliptic beta integrals on the $A_n$ root system and build corresponding
Hamiltonians (all of which coincide in the rank 1 case). Although all these
models are degenerate---their particles' pairwise coupling constant
is fixed in one or another way, it would be interesting to
clarify whether these models define new completely integrable
quantum systems.

\bigskip
\centerline{A\scriptsize CKNOWLEDGMENTS}
\medskip
These lectures are dedicated to the memory of a prominent russian
mathematician A. A. Bolibrukh whom I knew for a decade and
who has untimely passed away one year before this workshop.
He considered elliptic beta integrals as some universal
mathematical objects and presented the paper
\cite{spi:beta} for publication as a communication of the
Moscow mathematical society.

My sincere thanks go to the organizers of this workshop M. Noumi,  K. Takasaki
and the RIMS staff for the invitation to lecture at it and kind hospitality.
I am deeply indebted also to M. Ito for preparation of the
first draft of these lecture notes.

{\footnotesize

}
\end{document}